\documentclass[11pt]{amsart}
\usepackage{fullpage}
\usepackage{amsmath,amssymb,amsthm,graphicx,url}
\usepackage{booktabs,tabularx}
\usepackage{epsfig}

\theoremstyle{definition}
\newtheorem{theorem}{Theorem}
\newtheorem{lemma}[theorem]{Lemma}
\newtheorem{cor}[theorem]{Corollary}
\newtheorem{conj}[theorem]{Conjecture}
\newtheorem{prop}[theorem]{Proposition}
\newtheorem{example}[theorem]{Example}
\newtheorem{definition}[theorem]{Definition}

\newtheorem{remark}[theorem]{Remark}
\newtheorem{question}[theorem]{Question}

\DeclareMathOperator{\Par}{Par}

\DeclareMathOperator{\wt}{wt}
\DeclareMathOperator{\drw}{drw}

\newcommand{\dqbin}[3]{\displaystyle\genfrac{[}{]}{0pt}{}{#1}{#2}_{#3}}
\newcommand{\tqbin}[3]{\textstyle\genfrac{[}{]}{0pt}{}{#1}{#2}_{#3}} 
\newcommand{\CC}{\mathbb{C}}

\newcommand{\ZZ}{\mathbb{Z}}
\newcommand{\NN}{\mathbb{N}}

\newcommand{\RRR}{\mathsf{R}}

\newcommand{\rr}{b} 
\newcommand{\s}{a} 

\newcommand{\DR}{\mathsf{DR}}
\newcommand{\PF}{\mathsf{PF}}

\newcommand{\Cat}{\mathsf{Cat}}
\newcommand{\x}{\mathbf{x}}
\newcommand{\y}{\mathbf{y}}

\newcommand{\bounce}{\mathsf{bounce}}

\newcommand{\area}{\mathsf{area}}
\newcommand{\areap}{\mathsf{area'}}
\newcommand{\dinv}{\mathsf{dinv}}

\newcommand{\IDes}{\mathsf{IDes}}
\newcommand{\ml}{\mathsf{ml}}
\newcommand{\QSym}{\mathsf{QSym}}

\newcommand{\Sym}{\mathsf{Sym}}
\newcommand{\zetamap}{\mathsf{zeta}}

\DeclareMathOperator{\maj}{maj}

\newcommand{\arm}{\mathsf{arm}}
\newcommand{\leg}{\mathsf{leg}} 
\newcommand{\Hilb}{\mathsf{Hilb}}
\newcommand{\Frob}{\mathsf{Frob}}
\newcommand{\sweep}{\mathsf{sweep}}

\newcommand{\Schro}{{\sf Sch}} 
\newcommand{\patR}{\mathcal{R}} 
\newcommand{\patD}{\mathcal{D}} 
\newcommand{\ptnR}{\mathrm{R}}  
\newcommand{\ptnD}{\mathrm{D}}  

\makeatletter
\def\imod#1{\allowbreak\mkern10mu({\operator@font mod}\,\,#1)}
\makeatother

\begin{document}

\title{Rational Parking Functions and Catalan Numbers}

\date{\today}

\author{Drew Armstrong}
\address{Dept. of Mathematics\\
  University of Miami \\
  Coral Gables, FL 33146}
\email{armstrong@math.miami.edu}

\author{Nicholas A. Loehr}
\address{Dept. of Mathematics\\
  Virginia Tech \\
  Blacksburg, VA 24061-0123 \\
and Mathematics Dept. \\
United States Naval Academy \\
Annapolis, MD 21402-5002}
\email{nloehr@vt.edu}

\thanks{This work was partially supported by a grant from the Simons
 Foundation (\#244398 to Nicholas Loehr).}

\author{Gregory S. Warrington}
\address{Dept. of Mathematics and Statistics\\
  University of Vermont \\
  Burlington, VT 05401}
\email{gregory.warrington@uvm.edu}

\thanks{Third author supported in part by National Science Foundation
  grant DMS-1201312.}

\vspace*{.3in}

\begin{abstract} 
  The ``classical'' parking functions, counted by the Cayley number
  $(n+1)^{n-1}$, carry a natural permutation representation of the
  symmetric group $S_n$ in which the number of orbits is the Catalan
  number $\frac{1}{n+1}\binom{2n}{n}$.  In this paper, we will
  generalize this setup to ``rational'' parking functions indexed by a
  pair $(a,b)$ of coprime positive integers.  We show that these
  parking functions, which are counted by $b^{a-1}$, carry a
  permutation representation of $S_a$ in which the number of orbits is
  the ``rational'' Catalan number $\frac{1}{a+b}\binom{a+b}{a}$.  We
  compute the Frobenius characteristic of the $S_a$-module of
  $(a,b)$-parking functions.  Next we study $q$-analogues of the
  rational Catalan numbers, proposing a combinatorial formula for
  $\frac{1}{[a+b]_q}\tqbin{a+b}{a}{q}$ and relating this formula to a
  new combinatorial model for $q$-binomial coefficients. Finally, we
  discuss $q,t$-analogues of rational Catalan numbers and parking
  functions (generalizing the shuffle conjecture for the classical
  case) and present several conjectures.
\end{abstract}

\maketitle

\section{Introduction}
\label{sec:intro}


\subsection{Overview of Parking Functions}
\label{subsec:overview}

The goal of this paper is to generalize the theory of
``classical'' parking functions --- counted by $(n+1)^{n-1}$ --- to
the theory of ``rational'' parking functions --- counted by $b^{a-1}$
for coprime positive integers $a,b\in\mathbb{N}$. The combinatorial
foundation for this theory is given by rational Dyck paths, which are
lattice paths staying weakly above a line of slope $a/b$. The
classical case corresponds to $(a,b)=(n,n+1)$, which for most purposes
is equivalent to considering a line of slope $1$.

One major algebraic motivation for studying parking functions comes
from representation theory. 
It was conjectured by Haiman~\cite[Conjecture
2.1.1]{haiman-conjectures} that the $S_n$-module of ``diagonal
coinvariants'' has dimension $(n+1)^{n-1}$ and that the naturally
bi-graded character of this module might be encoded by parking
functions~\cite[Conjecture 2.6.3]{haiman-conjectures}.  The Shuffle
Conjecture~\cite[Conjecture 3.1.2]{HHLRU} (which is still open) gives
such a precise description.  We suggest in Definition~\ref{def:abpf} a
way to generalize the {\bf combinatorial} form of the Shuffle
Conjecture to rational parking functions (see also
Conjecture~\ref{conj:abpf}). However, we do not know what {\bf
  algebraic} or {\bf geometric} objects might underlie the generalized
combinatorics. The fact that the combinatorics works out so nicely
suggests that there must be an underlying reason.  
See~\cite{GMV,gorsky-negut,hikita} for related work that arises from a 
more geometric viewpoint.

The Shuffle Conjecture and related formulas for $q,t$-Catalan numbers
involve certain statistics on Dyck paths, namely $\area$ paired with either 
$\dinv$ or $\bounce$. The $\area$ statistic is natural and the $\dinv$ and 
$\bounce$ statistics are strange. However, there is a bijection $\zetamap$
on Dyck paths that sends $\dinv$ to $\area$ and $\area$ to $\bounce$. Thus one
could say that there is really just one natural statistic $\area$ and
one strange map $\zetamap$. In this paper we will define an analogous
map that we use to define $q,t$-analogues
of rational parking functions.  Experimentally, this map has beautiful
combinatorial properties.  Conjecture~\ref{conj:abpf} contains several
conjectures pertaining to the symmetry and $t=1/q$ specializations of
these rational $q,t$-parking functions. In a forthcoming paper~\cite{sweep}
we will show that our map belongs to a whole family of
\textbf{sweep maps} that contain a wealth of combinatorial information.

We will also study $q$-analogues of rational Catalan numbers that
are closely related to the $t=1/q$ specialization mentioned above.
These rational $q$-Catalan numbers are given algebraically by
$\frac{1}{[a+b]_q}\tqbin{a+b}{a,b}{q}$ where $\gcd(a,b)=1$
(see \S\ref{subsec:notation-qbin} for the definition of this notation).
Although these expressions have long been known to be polynomials
in $q$ with nonnegative integer coefficients, it is an open problem
to find a combinatorial interpretation for these polynomials.
We propose such a combinatorial model in \S\ref{sec:rat-qcat}
along with a related non-standard combinatorial interpretation
for general $q$-binomial coefficients.

\subsection{Outline of the Paper}
\label{subsec:outline}
The structure of the paper is as follows.  
We conclude the introduction by reviewing standard notation 
concerning partitions, symmetric functions, and $q$-binomial coefficients that 
will be used throughout the paper; for more details, 
consult~\cite{loehr-bij,Macd,stanley-ec2}.  \S\ref{sec:classical-pf}
gives background pertaining to diagonal coinvariants and
classical parking functions, deriving several formulas for the
ungraded Frobenius character of these modules.
\S\ref{sec:RPF} generalizes this discussion to the case of rational
parking functions. The key to deriving the Frobenius characters of
these new modules is Proposition~\ref{prop:multinomial}, which
enumerates rational parking functions with a specified vertical run
structure.  \S\ref{sec:rat-qcat} studies the rational $q$-Catalan
numbers, proposing new combinatorial formulas for these polynomials
and for general $q$-binomial coefficients based on certain partition
statistics. We prove that Conjecture~\ref{conj:nonstd-qbin}, which
proposes a novel combinatorial interpretation for $q$-binomial
coefficients, is equivalent to Conjecture~\ref{conj:rat-qcat}, which
proposes a combinatorial interpretation for rational $q$-Catalan
numbers.  \S\ref{sec:shuffle} reviews the classical theory of
$q,t$-parking functions and $q,t$-Catalan numbers, culminating in the
Shuffle Conjecture for the Frobenius series of the doubly-graded
module of diagonal coinvariants.  \S\ref{sec:qt-ratcat} generalizes
this theory to the case of rational $q,t$-Catalan numbers and rational
$q,t$-parking functions.  Finally, \S\ref{sec:comp} includes explicit
computations of various polynomials considered in this paper.

\subsection{Notation for Partitions}
\label{subsec:notation-ptn}

A \textbf{partition} is a weakly decreasing sequence
$\lambda=(\lambda_1,\lambda_2,\ldots,\lambda_s)$ of positive integers.
The $i$-th \textbf{part} of $\lambda$ is $\lambda_i$.  The
\textbf{area} of $\lambda$ is $|\lambda|=\lambda_1+\lambda_2+\cdots
+\lambda_s$. We say that $\lambda$ is a partition \textbf{of $n$}
(denoted by $\lambda\vdash n$) when $|\lambda|=n$.  The
\textbf{length} of $\lambda$ is $\ell(\lambda)=s$, the number of parts
of $\lambda$.  We consider the empty sequence to be the unique
partition of $0$; this partition has length zero.  It is sometimes
convenient to add one or more zero parts to the end of a partition;
this does not change the length of the partition.  Let $\Par$ be the
set of all partitions.  For all $j\geq 1$, $m_j(\lambda)$ is the
number of parts of $\lambda$ equal to $j$. The integer $z_{\lambda}$
is defined by
\[ z_{\lambda}=\prod_{j\geq 1} j^{m_j(\lambda)}m_j(\lambda)!. \]
For example, $\lambda=(4,4,1,1,1)$ is a partition with
$|\lambda|=11$, $\ell(\lambda)=5$, $m_1(\lambda)=3$, $m_4(\lambda)=2$,
and $z_{\lambda}=1^3 4^2 3!2!=192$.

\subsection{Notation for Symmetric Functions}
\label{subsec:notation-symm}

Let $\Sym$ denote the ring of symmetric functions, which we view as a
subring of $\CC[[x_1,x_2,x_3,\ldots ]]=\CC[[\x]]$ as
in~\cite[Ch. I]{Macd}. We now recall some bases of $\Sym$ used later
in the paper.  First, the \textbf{complete homogeneous symmetric
  functions} $h_i(\x)$ are defined by the generating function
\begin{equation}\label{eq:Ht}
H(t)=\sum_{i\geq 0} h_i(\x) t^i = \prod_{j\geq 1} \frac{1}{(1-x_jt)}.
\end{equation}
For any partition $\lambda$ of length $l$, let 
$h_{\lambda}(\x)=h_{\lambda_1}(\x)\cdots h_{\lambda_l}(\x)$.
Then $\{h_{\lambda}(\x):\lambda\in\Par\}$ is the
\textbf{complete homogeneous basis} of $\Sym$.  
Second, the \textbf{power sum symmetric functions} $p_i(\x)$
are defined for all $i\geq 1$ by setting 
$p_i(\x)=x_1^i+x_2^i+\cdots+x_k^i+\cdots$.
For any partition $\lambda$ of length $l$, let 
$p_{\lambda}(\x)=p_{\lambda_1}(\x)\cdots p_{\lambda_l}(\x)$.
Then $\{p_{\lambda}(\x):\lambda\in\Par\}$ is the
\textbf{power sum basis} of $\Sym$.  
Third, the \textbf{monomial symmetric function} $m_{\lambda}(\x)$
is the sum of all distinct monomials 
$x_{i_1}^{\lambda_1}x_{i_2}^{\lambda_2}\cdots x_{i_l}^{\lambda_l}$
with nonzero exponents $\lambda_1,\lambda_2,\ldots,\lambda_l$ 
(in some order). The \textbf{monomial basis} of $\Sym$ is
$\{m_{\lambda}(\x):\lambda\in\Par\}$.
Finally, the \textbf{Schur basis} of $\Sym$ (defined in~\cite{Macd})
is denoted $\{s_{\lambda}(\x):\lambda\in\Par\}$.

\subsection{Notation for $q$-Binomial Coefficients}
\label{subsec:notation-qbin}

Let $q$ be a formal variable.  For all integers $n\geq k\geq 0$,
define the \textbf{$q$-integer} $[n]_q=1+q+q^2+\cdots+q^{n-1}$, the
\textbf{$q$-factorial} $[n]!_q=\prod_{j=1}^n [j]_q$, and the
\textbf{$q$-binomial coefficient}
\begin{equation}\label{eq:qtbin}
  \dqbin{n}{k}{q}=\dqbin{n}{k,n-k}{q}=\frac{[n]!_q}{[k]!_q[n-k]!_q}.
\end{equation}
We usually use the multinomial coefficient notation
$\tqbin{a+b}{a,b}{q}=\frac{[a+b]!_q}{[a]!_q[b]!_q}$ when discussing
$q$-binomial coefficients.

\section{Classical Parking Functions}
\label{sec:classical-pf}

The permutations of the set $[n]=\{1,2,\ldots,n\}$ are counted by the
factorial $n!$. Two other important objects in combinatorics,
trees and parking functions, are
counted by the Cayley number $(n+1)^{n-1}$. One can algebraically
motivate the progression from $n!$ to $(n+1)^{n-1}$ as follows.

\subsection{Diagonal Coinvariants}
\label{sec:diag-coinv}

The symmetric group $S_n$ acts on the polynomial ring $\CC[x_1,x_2,\ldots,x_n]$ 
by permuting variables. 
Isaac Newton knew that the subring of ``symmetric polynomials'' is
generated by the \textbf{power sum polynomials} 
\[ p_k=p_k(x_1,\ldots,x_n)=x_1^k+x_2^k+\cdots+x_n^k\qquad (1\leq k\leq n). \]
Claude Chevalley~\cite{chevalley} knew that the quotient ring of {\bf
  coinvariants}
\begin{equation*}
\RRR_n=\CC[x_1,x_2,\ldots,x_n]/\langle p_1,p_2,\ldots,p_n\rangle
\end{equation*}
is isomorphic to a graded version of the regular representation of $S_n$, 
and so has dimension $n!$. Moreover, the Hilbert series of $\RRR_n$ is 
the $q$-factorial:
\begin{equation*}
  \Hilb_{\RRR_n}(q)=\sum_{i\geq 0} \dim(\RRR_n^{(i)})\, q^i 
= \prod_{j=1}^n (1+q+\cdots +q^{j-1}) = [n]!_q.
\end{equation*}
Armand Borel~\cite{borel} knew that $\RRR_n$ is also the cohomology ring
of the complete flag variety. It turns out that the structure of
$\RRR_n$ is closely related to the combinatorial structure of
permutations.

More generally, $S_n$ acts diagonally on the polynomial ring
$\CC[x_1,x_2,\ldots,x_n,y_1,y_2,\ldots,y_n]$ by {\em simultaneously}
permuting the $x$-variables and the $y$-variables. Hermann
Weyl~\cite{weyl} knew that the subring of $S_n$-invariant polynomials
is generated by the \textbf{polarized power sums}:
\begin{equation*}
p_{k,\ell} = p_{k,\ell}(x_1,\ldots,x_n,y_1,\ldots,y_n)=
x_1^ky_1^\ell+x_2^ky_2^\ell+\cdots+x_n^ky_n^\ell
\qquad(k,\ell\in\NN,k+\ell>0).
\end{equation*}
Mark Haiman~\cite[Conjecture 2.1.1]{haiman-conjectures} conjectured
that the quotient ring of {\bf diagonal coinvariants}
\begin{equation*}
\DR_n= \CC[x_1,x_2,\ldots,x_n,y_1,y_2,\ldots,y_n]/\langle 
  p_{k,\ell}: k+\ell >0\rangle
\end{equation*}
has dimension $(n+1)^{n-1}$ as a vector space over $\CC$. Furthermore,
this vector space is bi-graded by $\x$-degree and $\y$-degree, and he
saw hints that the bi-graded Hilbert series
\begin{equation*}
\mathrm{Hilb}_{\DR_n}(q,t)=\sum_{i,j\geq 0} \dim(\DR_n^{(i,j)})\, q^i t^j
\end{equation*}
is closely related to well-known structures in combinatorics.  Haiman,
along with Adriano Garsia, made several conjectures in this 
direction~\cite{GH-qtcat}. However, because the polarized
power sums $p_{k,\ell}$ are not algebraically independent, it turned
out to be quite difficult to prove these conjectures. Some have now
been proved, some are still open, and in general the subject remains
very active.

\subsection{Parking Functions}
\label{subsec:PF}

Arthur Cayley~\cite{cayley} showed that there are $(n+1)^{n-1}$ trees with
$n+1$ labeled vertices. However, for the purposes of studying diagonal
coinvariants, we prefer to discuss \emph{parking 
functions}~\cite{konheim-weiss,riordan},  which can be defined as follows.
There are $n$ cars $C_1,C_2,\ldots, C_n$ that want to park in $n$ spaces 
along a one-way street. Each car $C_i$ has a preferred spot $a_i\in [n]$ 
and the cars park in order: $C_1$ first, $C_2$ second, etc. When car $C_i$ 
arrives it will try to park in spot $a_i$. If spot $a_i$ is already taken 
it will park in the first available spot after $a_i$. If no such spot exists, 
the parking process fails. We say that the $n$-tuple 
$(a_1,a_2,\ldots,a_n)\in [n]^n$ is a {\bf classical parking function} 
if it allows all of the cars to park. 

For example, $(a_1,a_2,a_3,a_4,a_5)=(2,4,1,2,1)$ is a parking 
function. The order of the parked cars is shown here:
\bigskip
\begin{center}
\includegraphics[scale=.6]{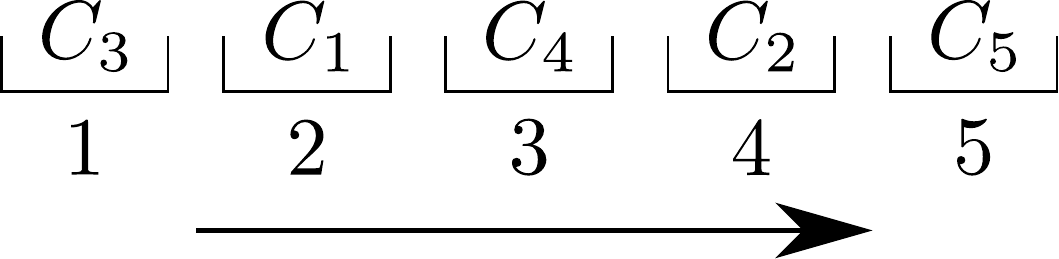}
\end{center}

One may check that $(a_1,a_2,\ldots,a_n)\in[n]^n$ is a parking
function if and only if its increasing rearrangement $b_1\leq b_2\leq
\cdots \leq b_n$ satisfies $b_i\leq i$ for all $i\in [n]$. This leads
to a few observations.

First, the set of parking functions $(a_1,a_2,\ldots,a_n)$
is closed under permuting subscripts.  Let $\PF_n$ denote both the set
of parking functions and the corresponding permutation representation
of $S_n$~\cite[Def. 1.3.2]{sagan}.  Haiman conjectured that, after
a sign twist, the ungraded $S_n$-module $\DR_n$ of diagonal coinvariants 
is isomorphic to $\PF_n$. We describe the graded version of this
conjecture in~\S\ref{subsec:graded-pfn} below.

Second, the {\em increasing} parking functions $(a_1\leq
a_2\leq \cdots \leq a_n)$ are in bijection with the set of Dyck
paths. We define a {\bf classical Dyck path of order $n$} as a sequence in
$\{N,E\}^n$ with the property that every initial subsequence has 
at least as many $N$'s as $E$'s. If we read from left to right,
interpreting $N$ as ``go north'' and $E$ as ``go east'', we can think
of this as a lattice path in $\ZZ^2$ from $(0,0)$ to $(n,n)$ staying
weakly above the diagonal line $y=x$. Now, given an increasing parking function
$P=(a_1\leq a_2\leq \cdots \leq a_n)$, let $r_i$ be the number of
times that $i$ occurs in $P$. Then we associate $P$ with the Dyck path
\begin{equation*}
\overbrace{N\cdots N}^{r_1}E
\overbrace{N\cdots N}^{r_2}E\cdots E
\overbrace{N\cdots N}^{r_n}E.
\end{equation*}
For example, the increasing parking function
$(1\leq1\leq2\leq2\leq4)\in\PF_5$ corresponds to the Dyck path
$NNENNEENEE$:
\begin{center}
\includegraphics[scale=.5]{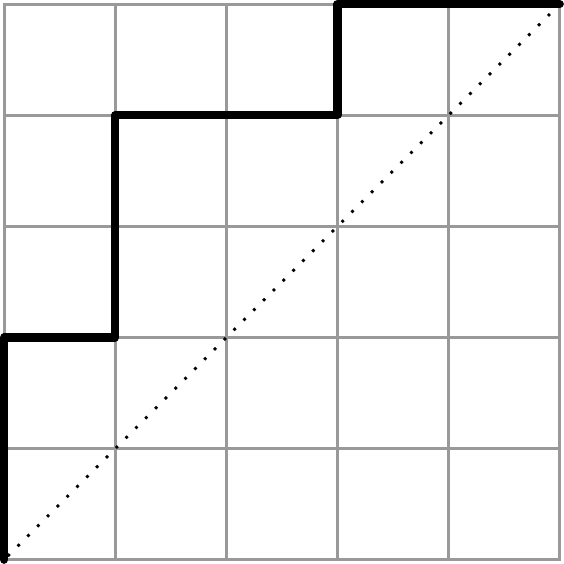}
\end{center}

The number of Dyck paths of order $n$ is the \textbf{Catalan number}
\begin{equation*}
\Cat_n=\frac{1}{n+1}\binom{2n}{n}.
\end{equation*}

Third, generalizing the previous observation, we associate (possibly
non-increasing) parking functions with {\bf labeled Dyck paths}. Given
a parking function $P=(a_1,a_2,\ldots,a_n)$ we first draw the Dyck
path corresponding to the increasing rearrangement. The $i$-th
vertical run of the path has length $r_i$ (possibly zero), which
corresponds to the number of occurrences of $i$ in $P$.  If we have
$r_i=k$ with $i=a_{j_1}=a_{j_2}=\cdots=a_{j_k}$ then we will label the
$i$-th vertical run by the set of indices $\{j_1,j_2,\ldots,j_k\}$.
We do this by filling the boxes to the right of the $i$-th vertical
run with the labels $j_1<j_2<\cdots <j_k$, increasing vertically. For
example, the parking function $(2,4,1,2,1)\in\PF_5$ corresponds to the
following labeled Dyck path: 
\begin{center}
\includegraphics[scale=.5]{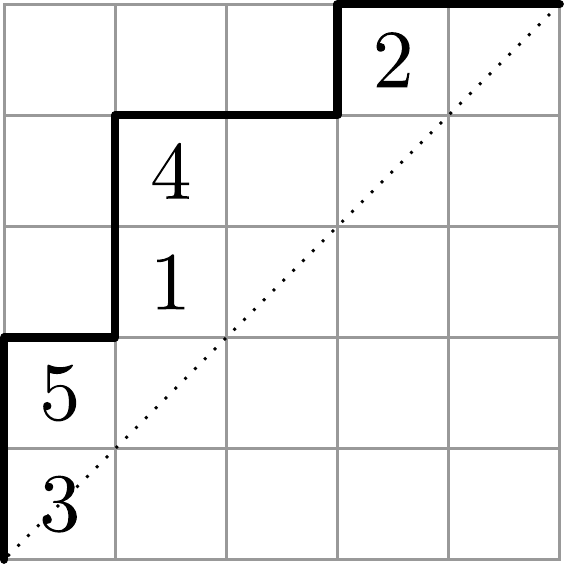}
\end{center}
This is the standard way to encode parking functions in the literature
on diagonal coinvariants. We will use this encoding to compute the
structure of the $S_n$-module $\PF_n$.

\subsection{Frobenius Characteristic of $\PF_n$}
\label{subsec:frob-pf}

We have observed that $S_n$ acts on parking functions $(a_1,\ldots,a_n)$
by permuting subscripts. Translating this to an action on labeled Dyck paths,
$S_n$ acts on a labeled path by permuting labels and then 
reordering labels in each column so that labels still increase
reading up each column. From this description, we see that the orbits
of this action
are in bijection with Dyck paths, hence are counted by the Catalan
number $\Cat_n=\frac{1}{n+1}\binom{2n}{n}$. 
Suppose $\lambda\vdash n$ and a given Dyck path has
$m_i=m_i(\lambda)$ vertical runs of length $i$ (for $1\leq i\leq n$).
Then the orbit corresponding to this Dyck path has
stabilizer isomorphic to the Young subgroup
$S^{\lambda}\leq S_n$, where
\begin{equation}\label{eq:young-subgp}
S^{\lambda}\cong S_1^{m_1}\times S_2^{m_2} \times \cdots \times S_n^{m_n}.
\end{equation}
The number of Dyck paths with this vertical run structure is known to equal
\begin{equation}\label{eq:count-by-runs}
\frac{1}{n+1}\binom{n+1}{m_0,m_1,m_2,\ldots,m_n} 
= \frac{n!}{m_0!m_1!m_2!\cdots m_n!},
\end{equation}
where we define $m_0$ so that $\sum_{i=0}^n m_i = n+1$; more
specifically, $m_0=n+1-\ell(\lambda)$.
(Equation~\ref{eq:count-by-runs} follows, for instance, from a more
general result proved in Proposition \ref{prop:multinomial} below.)

Next, recall that the {\bf Frobenius characteristic map} 
is an isomorphism from the $\CC$-algebra of symmetric group
representations to the $\CC$-algebra $\Sym$. This map
sends a class function $\chi:S_n\rightarrow\CC$ to
the symmetric function $\sum_{\mu\vdash n} 
\chi(\mu)\frac{p_{\mu}(\x)}{z_{\mu}}$, and the map sends
the irreducible character $\chi^{\lambda}$ to the
Schur function $s_{\lambda}(\x)$~\cite[\S4.7]{sagan}.
Let $\Frob_{\PF_n}$ denote the image of the character of $\PF_n$
under the Frobenius characteristic map.

Each orbit of $\PF_n$ with stabilizer isomorphic to $S^{\lambda}$ contributes 
a term $h_{\lambda}(\x)$ to $\Frob_{\PF_n}$~\cite[p. 113--114]{Macd}.
The number of such orbits is given by~\eqref{eq:count-by-runs}, and hence 
\begin{align}\label{eq:frob-h}
  \Frob_{\PF_n} &= \sum_{\lambda\vdash n} 
  \frac{1}{n+1}\binom{n+1}{n-\ell(\lambda)+1,m_1(\lambda),\ldots,m_n(\lambda)}
  h_{\lambda}(\x)\\
  &= \frac{1}{n+1}\sum_{\lambda\vdash n} 
  \binom{n+1}{n-\ell(\lambda)+1,m_1(\lambda),\ldots,m_n(\lambda)} 
  h_1(\x)^{m_1(\lambda)}\cdots h_n(\x)^{m_n(\lambda)}.
\end{align}
Using the multinomial theorem, the generating function~\eqref{eq:Ht},
and the fact that $h_0(\x)=1$, we obtain the identity
\begin{equation}\label{eq:frob1}
  \Frob_{\PF_n} = \left.\frac{1}{n+1} [H(t)]^{n+1}\right|_{t^n}.
\end{equation}

\subsection{Character of $\PF_n$.}
\label{subsec:char-pfn}

To extract information from~\eqref{eq:frob1}, we recall the 
{\bf Cauchy product}
\begin{equation}\label{eq:pixy}
  \Pi(\x,\y)= \prod_{i=1}^\infty \prod_{j=1}^\infty \frac{1}{(1-x_iy_j)}.
\end{equation}
The following well-known lemma shows how the Cauchy product
can be used to detect dual bases of $\Sym$.

\begin{lemma}~\cite[I.4.6]{Macd},\cite[Thm. 10.131]{loehr-bij}
\label{lem:ip} 
Define the \textbf{Hall inner product} $\langle\cdot,\cdot\rangle$ on $\Sym$ 
by requiring that $\{h_\lambda:\lambda\in\Par\}$ 
and $\{m_\lambda:\lambda\in\Par\}$ be dual bases. 
For any two bases $\{v_{\lambda}:\lambda\in\Par\}$
and $\{w_{\lambda}:\lambda\in\Par\}$ of $\Sym$
such that all $v_{\lambda}$ and $w_{\lambda}$ are
homogeneous of degree $|\lambda|$, we have
  $\Pi(\x,\y) = \sum_{\lambda\in\Par} v_\lambda(\x) w_\lambda(\y)$
if and only if these two bases are dual with respect to the Hall inner product. 
In particular,
\begin{equation}\label{eq:newpixy}
 \Pi(\x,\y)=\sum_{\lambda\in\Par} h_{\lambda}(\x)m_{\lambda}(\y)
=\sum_{\lambda\in\Par} \frac{p_{\lambda}(\x)p_{\lambda}(\y)}{z_{\lambda}}
=\sum_{\lambda\in\Par} s_{\lambda}(\x)s_{\lambda}(\y). 
\end{equation}
\end{lemma}

It follows from~\eqref{eq:Ht} and~\eqref{eq:pixy} that by
setting $n+1$ of the $y$ variables equal to $t$ and the rest equal to
zero, $\Pi(\x,\y)$ specializes to $[H(t)]^{n+1}$.  Furthermore,
one sees directly that $p_{\lambda}(\y)$ specializes to
$(n+1)^{\ell(\lambda)} t^{|\lambda|}$.  
These facts, combined with~\eqref{eq:frob1} and~\eqref{eq:newpixy}, yield 
\begin{equation*}
\Frob_{\PF_n} = \left.\frac{1}{n+1}[H(t)]^{n+1}\right|_{t^n} 
= \left.\frac{1}{n+1}\sum_{\lambda\in\Par} (n+1)^{\ell(\lambda)} t^{|\lambda|} 
\frac{p_\lambda(\x)}{z_\lambda}\right|_{t^n}.
\end{equation*}
Hence
\begin{equation*}
  \Frob_{\PF_n} = \sum_{\lambda\vdash n} (n+1)^{\ell(\lambda)-1} \frac{p_\lambda(\x)}{z_\lambda}.
\end{equation*}
Applying the inverse of the Frobenius map, this formula tells
us the character $\chi$ of the $S_n$-module $\PF_n$. 
On one hand, for $w\in S_n$, $\chi(w)$ is the number of parking
functions in $\PF_n$ fixed by the action of $w$.
On the other hand, if $w$ has cycle type $\lambda$
(so that $w$ is a product of $\ell(\lambda)$ disjoint cycles),
the preceding formula shows that $\chi(w)= 
(n+1)^{\ell(\lambda)-1}$. In particular, taking $w$ to be the identity
permutation in $S_n$, which has cycle type $\lambda=(1^n)$ and fixes all 
parking functions in $\PF_n$, we see that $|\PF_n|=(n+1)^{n-1}$.
(There are easier ways to obtain this result, but
we wanted to illustrate the power of the generating function method.)

\subsection{Schur Expansion of $\PF_n$.}
\label{subsec:schur-pfn}

The derivation of the Frobenius series of $\PF_n$ can be redone
using the Schur symmetric functions instead of the power sum basis.
As above, we specialize $\y$ to consist of zeros along with $n+1$ copies 
of $t$. Since $s_{\lambda}(\y)$ is homogeneous of degree $|\lambda|$,
this specialization changes $s_{\lambda}(\y)$ into 
$t^{|\lambda|}s_{\lambda}(1^{n+1})$, where
$s_\lambda(1^{n+1})$ indicates that $n+1$ variables have been set to
$1$ and the rest to zero. By~\eqref{eq:newpixy},
\begin{equation*}
  \Frob_{\PF_n} = \sum_{\lambda\vdash n} 
  \frac{s_\lambda(1^{n+1})}{n+1} s_\lambda(\x).
\end{equation*}
This gives the decomposition of the parking function module into
irreducible constituents. In the special case of hook shapes
$\lambda=(k,1^{n-k})$ we can compute the coefficients explicitly (see
~\cite[p. 364]{stanley-ec2}):
\begin{equation*}
\frac{s_{(k,1^{n-k})}(1^{n+1})}{n+1} = 
 \frac{1}{n+1} \binom{n-1}{k-1} \binom{n+k}{n}.
\end{equation*}
These integers are called the {\bf Schr\"oder numbers}, and they also describe
the $f$-vector of the associahedron (see
Pak-Postnikov~\cite[formula~(1--6)]{pak-postnikov}). Setting $k=n$, we
verify that the multiplicity of the trivial representation in $\PF_n$
is the classical Catalan number:
\begin{equation*}
\frac{s_{(n)}(1^{n+1})}{n+1}=\Cat_n=\frac{1}{n+1}\binom{2n}{n}.
\end{equation*}
These computations were originally done by Haiman \cite{haiman-conjectures},
Pak-Postnikov~\cite{pak-postnikov}, and Stanley~\cite{stanley-parking}.
 
\subsection{The Graded Version of $\PF_n$.}
\label{subsec:graded-pfn} 

As mentioned earlier, it is known that the parking function module
$\PF_n$ and the diagonal coinvariant ring $\DR_n$ are isomorphic as
(ungraded) $S_n$-modules.  However, $\DR_n$ also comes with a
symmetric bi-grading by $\x$-degree and $\y$-degree, and it is natural
to ask whether we can explain this bi-grading in terms of $\PF_n$.
This is the content of the famous Shuffle Conjecture, described in
\S\ref{sec:shuffle} below. This section studies the simpler case where
we grade $\DR_n$ by the $\x$-degree only.

There is a straightforward combinatorial construction that turns the 
module $\PF_n$ into a singly-graded vector space.  The grading is defined
on basis vectors by the statistic $\area:\PF_n\to\NN$, which
is the number of boxes fully contained between the labeled Dyck path 
and the diagonal $y=x$ . For example, the following figure shows
that $\area(2,4,1,2,1)=5$: \bigskip
\begin{center}
\includegraphics[scale=.5]{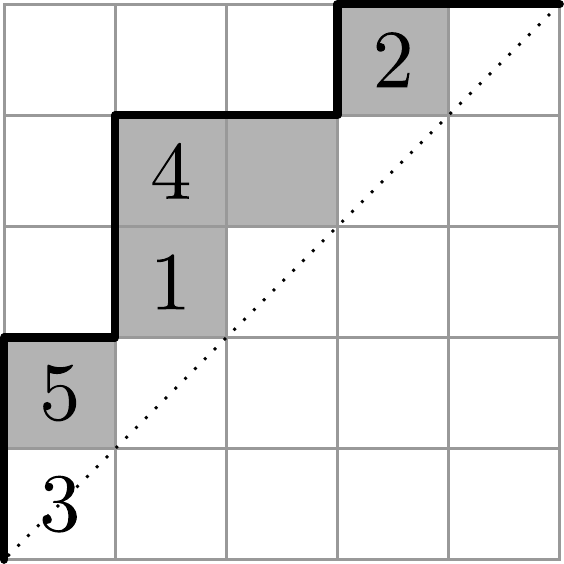}
\end{center}
Since permuting the labels of a labeled Dyck path leaves the area
unchanged, we have given $\PF_n$ the structure of a graded $S_n$-module.
Haiman conjectured and eventually
proved~\cite{GH-qtcat,haiman-conjectures,haiman-vanish} that this
module is isomorphic to $\epsilon\otimes\DR_n$ graded by $\x$-degree
only, where $\epsilon$ denotes the sign character.

\section{Rational Parking Functions}
\label{sec:RPF}

In this section we will generalize the $S_n$-module of \emph{classical}
parking functions $\PF_n$ labeled by a single positive integer
$n\in\NN$ to an $S_a$-module of \emph{rational} parking functions, denoted
$\PF_{a,b}$, labeled by two coprime positive integers $a$ and $b$. The
classical parking functions will correspond to the case
$(a,b)=(n,n+1)$. To define $\PF_{a,b}$ we must first discuss rational
Dyck paths.

\subsection{Rational Dyck Paths} 
\label{subsec:rat-dyck-paths}

Given positive integers $a,b\in\NN$,
let $\patR(N^a E^b)$ denote the subset of $\{N,E\}^{a+b}$ consisting
of words containing $a$ copies of $N$ and $b$ copies of $E$. We can
think of such a word as a lattice path in $\ZZ^2$ from $(0,0)$ to
$(b,a)$ by reading from left to right, interpreting $N$ as ``go
north'' and $E$ as ``go east.'' An element of $\patR(N^a E^b)$ is
called an {\bf $(a,b)$-Dyck path} if it stays weakly above the
diagonal line $y=\frac{a}{b}x$. Let $\patD(N^a E^b)$ denote the set of
$(a,b)$-Dyck paths. An example of a $(5,8)$-Dyck path is shown here:
\begin{center}
\includegraphics[scale=.5]{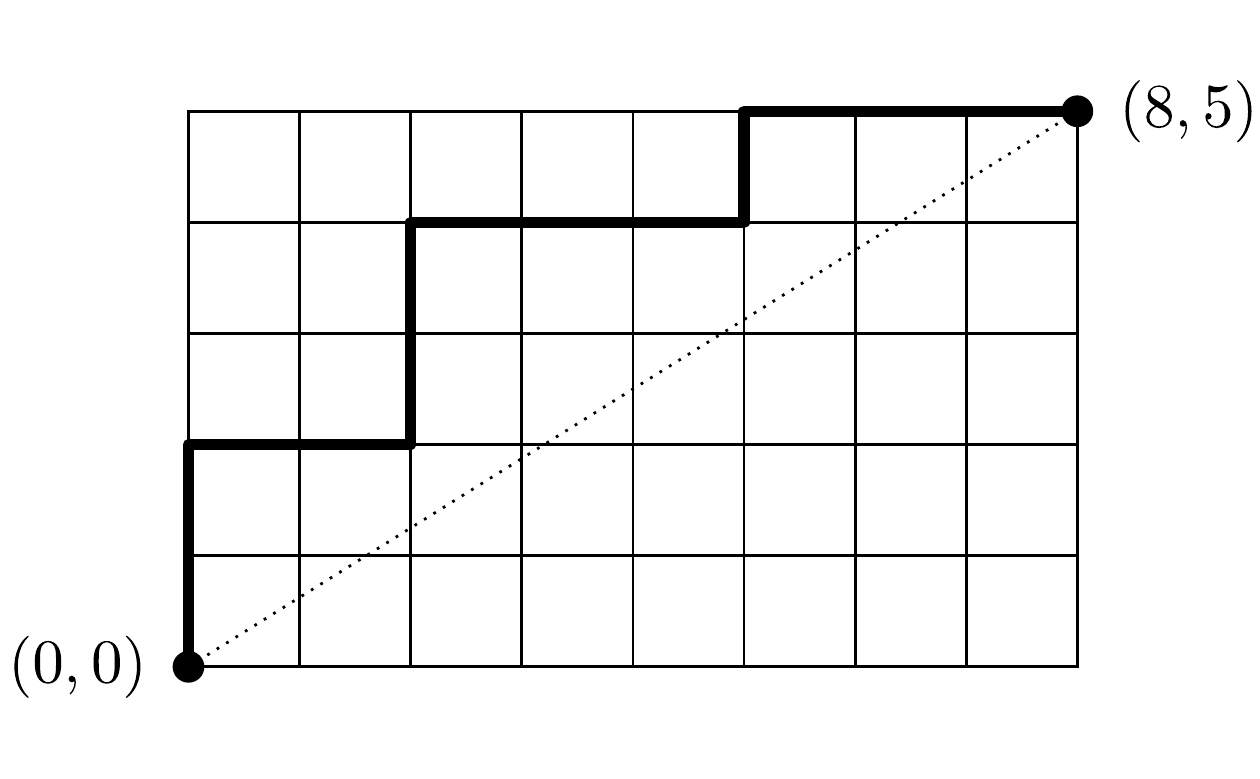}
\end{center}
Recall that the $(n,n)$-Dyck paths are called ``classical'' and they
are counted by the Catalan number
$\Cat_n=\frac{1}{n+1}\binom{2n}{n}$. Every
$(n,n+1)$-Dyck path must end with an east step, and by removing this east
step we obtain a canonical bijection
\begin{equation*}
\patD(N^n E^{n+1})\rightarrow \patD(N^nE^n).
\end{equation*}
In the case that $a,b$ are coprime, the enumerative theory of
Dyck paths is quite nice. In particular, it was known as early as 1954
(see Bizley~\cite{bizley}) that when $a,b$ are coprime we have
\begin{equation*}
  |\patD(N^a E^b)| 
  = \frac{1}{a+b}\binom{a+b}{a,\, b}=\frac{(a+b-1)!}{a!\, b!}.
\end{equation*}
We call this number the {\bf rational Catalan number}, 
denoted $\Cat_{a,b}$.  When $a,b$ are coprime, 
we call $\patD(N^a E^b)$ the set of {\bf rational Dyck
paths} in $\patR(N^a E^b)$. Observe that $\Cat_{n,n+1}=\Cat_n$.

The main result of this section is a formula enumerating
rational Dyck paths with a specified vertical run structure.
In this setting, a ``vertical run of length $i\geq 0$'' is
a subword $N^iE$ such that this subword is either preceded 
by $E$ or is the beginning of the entire word.
The classical (i.e., non-rational) version of this formula goes back
at least to Kreweras' work in 1972~\cite[Theorem 4]{kreweras} in the
context of noncrossing partitions.

\begin{prop}
\label{prop:multinomial} 
Let $a$ and $b$ be coprime positive integers.
Let $m_0,m_1,m_2,\ldots,m_a$ be nonnegative integers such that
$\sum_{i\geq 0} im_i=a$ and $\sum_{i\geq 0} m_i=b$.  Then the number of
$(a,b)$-Dyck paths in $\patD(N^aE^b)$ with $m_i$
vertical runs of length $i$ is given by
\[\frac{1}{b}\binom{b}{m_0,m_1,\ldots,m_a}
 =\frac{(b-1)!}{m_0!m_1!\cdots m_a!}.\] 
\end{prop}
\begin{proof}
  \emph{Step 1.} Let $Y$ be the set of paths $\pi\in\patR(N^aE^b)$
  such that $\pi$ has $m_i$ vertical runs of length $i$ for all $i\geq
  0$, and $\pi$ ends in an east step.  
  We show $$|Y|=\binom{b}{m_0,m_1,\ldots,m_a}.$$ Let $X$ denote the
  set $\patR(v_0^{m_0}v_1^{m_1}\cdots v_a^{m_a})$ of words containing
  $m_i$ copies of $v_i$ for each $i$.  Define
  $f:X\rightarrow\{N,E\}^*$ by replacing each letter $v_i$ in a word
  by $N^iE$. Since $\sum_i im_i=a$ and $\sum_i m_i=b$, we see that
  $f$ maps $X$ into $\patR(N^aE^b)$ and that $f(w)$ ends in an east
  step for all $w\in X$.  It is now clear that $f$ is a bijection of
  $X$ onto $Y$.  Since $|X|=\binom{b}{m_0,m_1,\ldots,m_a}$, Step 1 is
  complete.

\emph{Step 2.} Let $\pi \in Y$.  We associate a {\bf level} $l_i$
to the $i$-th lattice point on a path $\pi$ as follows.  Set $l_0 = 0$.  
For each $i > 0$, set 
\begin{equation*}
  l_i = 
  \begin{cases}
    l_{i-1} + b, & \text{ if the $i$-th step of $\pi$ is a north step};\\
    l_{i-1} - a, & \text{ if the $i$-th step of $\pi$ is an east step}.
  \end{cases}
\end{equation*}
Note that the level of a lattice point $(x,y)$ is $by-ax$.
We show that the levels $l_1,l_2,\ldots,l_{a+b}$ for a given $\pi\in Y$
are all \emph{distinct}.  For suppose $l_i=l_j$ for some $i<j$; we
prove that $i$ must be $0$ and $j$ must be $a+b$. Let $(x_i,y_i)$ and
$(x_j,y_j)$ be the lattice points reached by the $i$-th and $j$-th
steps of $\pi$.  We know $0\leq x_i\leq x_j\leq b$ and $0\leq y_i\leq
y_j\leq a$; also, $l_i=by_i-ax_i$ and $l_j=by_j-ax_j$. Since
$l_i=l_j$, $b(y_j-y_i)=a(x_j-x_i)$. But $a$ and $b$ are coprime, so
their least common multiple is $ab$. This forces $a$ to divide
$y_i-y_j$ and $b$ to divide $x_i-x_j$. So we must have $x_i=y_i=0$,
$x_j=b$, and $y_j=a$, giving $i=0$ and $j=a+b$ as claimed.

\emph{Step 3.} Define an equivalence relation $\sim$ on $Y$ by letting
$\pi_1\sim\pi_2$ iff there exist $w_1,w_2\in X$ such that
$f(w_1)=\pi_1$, $f(w_2)=\pi_2$, and $w_2$ is a cyclic shift of
$w_1$. We show that every equivalence class has size $b$ and contains
exactly one $(a,b)$-Dyck path. Fix $\pi\in Y$ and write $\pi=f(w)$ for
some $w\in X$. By cyclically shifting $w$ by $0,1,2,\ldots,b-1$, we
obtain a list of $b$ paths $\pi=\pi_0,\pi_1,\pi_2,\ldots,\pi_{b-1}$
(not yet known to be distinct), which are the paths in $Y$ equivalent
to $\pi$.  Suppose the list of steps in $\pi$ has east steps at
positions $i_1<i_2<\cdots<i_b=a+b$. Define $i_0=0$. Cyclically
shifting $w$ by $k$ steps, where $0\leq k<b$, has the effect of
cyclically shifting $\pi$ by $i_k$ steps. This cyclic shift will
replace the sequence of levels $l_0,l_1,\ldots$ for $\pi$ by the
new sequence
\begin{equation}\label{eq:shifted-wand-seq}
l_{i_k}-l_{i_k},l_{(i_k+1)}-l_{i_k},\ldots,
l_{a+b}-l_{i_k},l_1-l_{i_k},l_2-l_{i_k},\ldots,l_{i_k-1}-l_{i_k}.
\end{equation} 
If $m_0$ is the minimum level in $\pi$, it follows that $m_0-l_{i_k}$
is the minimum level in $\pi_k$. Since $l_{i_1},\ldots,l_{i_b}$ are
distinct by Step 2, it follows that the $b$ paths $\pi_0,\ldots,\pi_{b-1}$ 
all have distinct minimum levels, hence these $b$ paths must be pairwise
distinct. Furthermore, the minimum level $m_0$ in $\pi$ must occur at
the end of an east step, so $m_0=l_{i_j}$ for some $j$.  For every
$k$, the minimum level in~\eqref{eq:shifted-wand-seq} is
$l_{i_j}-l_{i_k}$, which is nonnegative iff $l_{i_j}\geq l_{i_k}$.
But since $l_{i_j}$ is the minimum level in $\pi$ and all levels are
distinct, the minimum level in $\pi_k$ is nonnegative iff $k=j$.  This
means that exactly one of the paths in the equivalence class of $\pi$
is an $(a,b)$-Dyck path.

\emph{Step 4.}
Step 3 shows that $Y$ decomposes into a disjoint union of $b$-element
subsets, each of which meets $\patD(N^aE^b)$ in exactly one
point. So the cardinality $|\patD(N^aE^b)|$ can be found by
dividing the multinomial coefficient in Step 1 by $b$.
\end{proof}

\subsection{Rational Parking Functions}
\label{subsec:rat-pf}

An {\bf $(a,b)$-parking function} is an $(a,b)$-Dyck path together
with a labeling of the north steps by the set $\{1,2,\ldots,a\}$ such
that labels increase in each column going north. For example, here is
a $(5,8)$-parking function.
\begin{center}
  \includegraphics[scale=1.0]{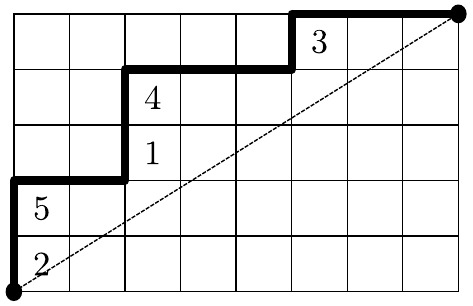}
\end{center}
The symmetric group $S_a$ acts on $(a,b)$-parking functions by
permuting labels and then reordering the labels within columns if
necessary.  Let $\PF_{a,b}$ denote the set, and also the $S_a$-module,
of $(a,b)$-parking functions.  Consistent with our terminology for
$(a,b)$-Dyck paths, we call $\PF_{a,b}$ the set of {\bf rational
  parking functions} in the case where $a$ and $b$ are coprime.  We make this
assumption of coprimality for the rest of this section.  In this case,
we can compute the Frobenius characteristic of $\PF_{a,b}$ by
the same method used in \S\ref{subsec:frob-pf} to compute the
Frobenius characteristic of $\PF_n$ (which is essentially $\PF_{n,n+1}$).

\begin{theorem} 
  The Frobenius characteristic of the $S_a$-module of parking
  functions $\PF_{a,b}$ has the following expansions in terms of
  complete homogeneous, power sum, and Schur symmetric functions, respectively:
  \begin{align*}
    \mathsf{Frob}(\PF_{a,b}) &= \sum_{\lambda\vdash a} \frac{1}{b}
    \binom{b}{b-\ell(\lambda),m_1(\lambda),m_2(\lambda),\ldots,m_a(\lambda)}
     h_{\lambda}(\x)
\\ &=\sum_{\lambda\vdash a} b^{\ell(\lambda)-1}\frac{p_\lambda(\x)}{z_\lambda}
\\ &=\sum_{\lambda\vdash a} \frac{s_\lambda(1^b)}{b} s_\lambda(\x).
  \end{align*}
  \label{th:rationalFrobenius}
\end{theorem}

\begin{proof}
The $h$-expansion follows from Proposition~\ref{prop:multinomial}
in the same way that~\eqref{eq:frob-h} followed from~\eqref{eq:count-by-runs}. 
(Note that the hypothesis that $a$ and $b$ are coprime is necessary here.) 
By the multinomial theorem, $\Frob_{\PF_{a,b}}$ is the coefficient of $t^a$ in 
the generating function $\frac{1}{b}[H(t)]^b$.  
The power sum and Schur expansions then follow 
from~\eqref{eq:pixy} and~\eqref{eq:newpixy}, 
by the same arguments used in
\S\ref{subsec:char-pfn} and \S\ref{subsec:schur-pfn}.  
\end{proof}

Recall that the $p$-expansion tells us the character of the
permutation action of $S_a$ on $\PF_{a,b}$. We can express this as
follows.

\begin{cor}
Let $w\in S_a$ be a permutation with $k$ cycles. Then the number of
elements of $\PF_{a,b}$ fixed by $w$ equals $b^{k-1}$. In
particular, taking $w$ to be the identity permutation in $S_a$
(which has $a$ cycles), $|\PF_{a,b}|=b^{a-1}$.  
\end{cor}

Using~\cite[page 364]{stanley-ec2}, we obtain the following 
formula for the coefficient of $s_{\lambda}$ when $\lambda$
is a hook shape.  

\begin{cor}
\label{cor:hook}
For $0\leq k\leq a-1$, the multiplicity of the hook Schur function 
$s_{(k+1,1^{a-k-1})}(\x)$ in $\Frob_{\PF_{a,b}}$ is
\begin{equation*}
\frac{s_{(k+1,1^{a-k-1})}(1^b)}{b} = \frac{1}{b} \binom{a-1}{k}\binom{b+k}{a}.
\end{equation*}
We call these integers the {\bf rational Schr\"oder numbers},
denoted $\Schro_{a,b;k}$. We note that
these numbers also occur as the $f$-vector of the recently discovered
\emph{rational associahedron}~\cite{ARW}.
\end{cor}

We already knew a special case of this result. Namely, when $k=a-1$, we
find that the multiplicity of the trivial character in $\PF_{a,b}$ is
the rational Catalan number
\begin{equation*}
{\sf Cat}_{a,b} = \frac{1}{b}\binom{b+a-1}{a,b-1}=\frac{(a+b-1)!}{a!\,b!}.
\end{equation*}
But since $\PF_{a,b}$ is a permutation module, the multiplicity of the
trivial character is the number of orbits. Since the orbits are
represented by $(a,b)$-Dyck paths, we recover Bizley's result
regarding the number of $(a,b)$-Dyck paths~\cite{bizley}.
This result can also be proved by an argument similar to the one
in Proposition~\ref{prop:multinomial}; see~\cite[\S 12.1]{loehr-bij}.

Observe from Corollary~\ref{cor:hook} that the sign
character of $S_a$ occurs in $\PF_{a,b}$ if and only if $b\geq a$.
More precisely, we see that the smallest value of $k$ for which
$s_{(k+1,1^{a-k-1})}(\x)$ occurs in $\Frob_{\PF_{a,b}}$ is $k=\max\{0,a-b\}$.

\subsection{Graded Version of $\PF_{a,b}$.}
\label{subsec:PFab-graded}

It is straightforward to generalize the $\area$ statistic on
classical Dyck paths to rational Dyck paths. For coprime $a,b\in\NN$
and $D\in \patD(N^aE^b)$, $\area(D)$ is the number of boxes fully
contained between $D$ and the diagonal line $y=\frac{a}{b}x$.
For example, the rational Dyck path displayed at the beginning
of~\S\ref{subsec:rat-pf} has area 5. The action of $S_a$ on $\PF_{a,b}$
preserves the $\area$ statistic, so this statistic turns $\PF_{a,b}$
into a singly-graded $S_a$-module, as in the classical case.

When $a$ and $b$ are coprime, the maximum value of $\area(D)$ over all
$(a,b)$-Dyck paths is $(a-1)(b-1)/2$. Indeed, the diagonal of the
$a\times b$ rectangle intersects a ribbon of $a+b-1$ boxes, as shown
here:
\begin{center}
\includegraphics[scale=1]{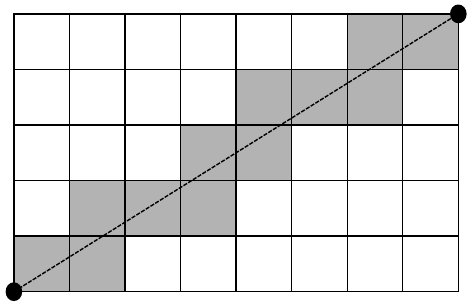}
\end{center}
Note that this ribbon divides the rest of the rectangle into two equal 
pieces of size
\begin{equation*}
\frac{ab-(a+b-1)}{2} = \frac{(a-1)(b-1)}{2}.
\end{equation*}

\section{Rational $q$-Catalan Numbers}
\label{sec:rat-qcat}

This section studies a $q$-analogue of the rational Catalan number
$\Cat_{a,b}$ obtained by replacing $\frac{1}{a+b}\binom{a+b}{a,b}$
by $\frac{1}{[a+b]_q}\tqbin{a+b}{a,b}{q}$. We conjecture combinatorial
interpretations for these polynomials, as well as a new combinatorial
formula for all $q$-binomial coefficients,
based on certain partition statistics. The main theorem of this section
shows that the conjecture for $q$-binomial coefficients implies the
conjecture for rational $q$-Catalan numbers.

\subsection{Background on $q$-Binomial Coefficients.}
\label{subsec:qbinom}

The $q$-binomial coefficients defined by~\eqref{eq:qtbin} are in fact 
polynomials in $\NN[q]$.  One can prove this algebraically by checking that
the $q$-binomial coefficients satisfy the recursions
\begin{equation}\label{eq:qbin-rec}
 \dqbin{\s+\rr}{\s,\rr}{q} = q^\s\dqbin{\s+\rr-1}{\s,\rr-1}{q} 
+\dqbin{\s+\rr-1}{\s-1,\rr}{q}
= \dqbin{\s+\rr-1}{\s,\rr-1}{q} +q^\rr\dqbin{\s+\rr-1}{\s-1,\rr}{q} 
\quad(\s,\rr>0)
\end{equation}
and initial conditions $\tqbin{\s}{\s,0}{q}=\tqbin{\rr}{0,\rr}{q}=1$. 

Recall the following
well-known combinatorial interpretation of $q$-binomial coefficients.
Let $\ptnR(\s,\rr)$ denote the set of 
partitions whose diagrams
fit in the box with corners $(0,0)$, $(\rr,0)$, $(\rr,\s)$, and
$(0,\s)$. More precisely, $\mu\in \ptnR(\s,\rr)$ iff
$\mu=(\mu_1,\mu_2,\ldots,\mu_\s)$ where the parts $\mu_j$ are integers
satisfying $\rr\geq \mu_1\geq \mu_2\geq\cdots\geq \mu_\s\geq 0$.  
Note that $|\mu|=\mu_1+\mu_2+\cdots+\mu_\s$ is the number of boxes in the
diagram of $\mu$.  
One may verify that
\begin{equation}\label{eq:qbin-std}
 \dqbin{\s+\rr}{\s,\rr}{q}=\sum_{\mu\in \ptnR(\s,\rr)} q^{|\mu|}
\end{equation}
by showing that the right side of~\eqref{eq:qbin-std}
satisfies the recursions in~\eqref{eq:qbin-rec} (see, 
e.g.,~\cite[\S6.7]{loehr-bij}). We refer to~\eqref{eq:qbin-std}
as the \emph{standard} combinatorial interpretation
for the $q$-binomial coefficients, to distinguish it
from the non-standard interpretation presented below.

\subsection{Rational $q$-Catalan Numbers}
\label{subsec:rat-qcat}

Suppose $\s$ and $\rr$ are positive integers with $\gcd(\s,\rr)=1$.  
The \emph{rational $q$-Catalan number for the slope $\s/\rr$} is
\begin{equation}\label{eq:qpoly}
 \Cat_{\s,\rr}(q)=\frac{1}{[\s+\rr]_q}\dqbin{\s+\rr}{\s,\rr}{q}
  =\frac{[\s+\rr-1]!_q}{[\s]!_q[\rr]!_q}
  =\frac{1}{[\rr]_q}\dqbin{\s+\rr-1}{\s,\rr-1}{q}
  =\frac{1}{[\s]_q}\dqbin{\s+\rr-1}{\s-1,\rr}{q}.
\end{equation}
It is known that $\Cat_{\s,\rr}(q)$ is a polynomial with coefficients
in $\NN$.  Haiman gives an algebraic proof of this fact 
in~\cite[Prop. 2.5.2]{haiman-conjectures} and also gives an
algebraic interpretation for these polynomials as the
Hilbert series of a suitable quotient ring of a polynomial ring
~\cite[Prop. 2.5.3 and 2.5.4]{haiman-conjectures}. 
This polynomial also appears to be connected to certain
modules arising in the theory of rational Cherednik algebras.
Our purpose here is to propose a combinatorial interpretation
for the rational $q$-Catalan numbers.  More precisely, given $\s,\rr>0$
with $\gcd(\s,\rr)=1$, our problem is to find a set $X(\s,\rr)$ of 
combinatorial 
objects and a statistic $\wt:X(\s,\rr)\rightarrow\NN$ such that
\[ \Cat_{\s,\rr}(q)=\sum_{w\in X(\s,\rr)} q^{\wt(w)}. \]

For certain special choices of $\s$ and $\rr$, this problem has been solved.
For instance, when $\s=n$ and $\rr=n+1$, we may take $X(\s,\rr)$ to be the
set $DW_n$ of \emph{Dyck words} consisting of $n$ zeroes and $n$ ones, such 
that every prefix of the word has at least as many zeroes as ones. 
(Dyck words are obtained from classical Dyck paths 
by replacing each $N$ by 0 and each $E$ by 1.) Given
a word $w=w_1w_2\cdots w_{2n}\in DW_n$, define the \emph{major index}
of $w$, denoted $\maj(w)$, to be the sum of all $i<2n$ with $w_i>w_{i+1}$.
MacMahon~\cite{macmahon} proved that
\[ \sum_{w\in DW_n} q^{\maj(w)}=\frac{1}{[n+1]_q}\dqbin{2n}{n,n}{q}
  =\Cat_{n,n+1}(q). \]
Another known case is $\s=n$ and $\rr=mn+1$, where $m,n$ are fixed positive
integers. In this case, we may take $X(\s,\rr)$ to be the set of lattice
paths from $(0,0)$ to $(mn,n)$ that never go below the line $x=my$.
The statistic needed here is obtained by taking the $t=1/q$ specialization
of the higher-order $q,t$-Catalan numbers~\cite{GH-qtcat}.
Using the $\area$ and $\bounce$ statistics defined in~\cite{loehr-mcat},
it is shown in~\cite[\S3.3]{loehr-mcat} that
\[ \sum_{w\in X(n,mn+1)} q^{mn(n-1)/2+\area(w)-\bounce(w)}
  =\frac{1}{[mn+1]_q}\dqbin{mn+n}{mn,n}{q}=\Cat_{n,mn+1}(q). \]
One can replace the power of $q$ here by $mn(n-1)/2+h(w)-\area(w)$,
where $h$ (also called $\dinv_m$) is defined in~\cite{loehr-mcat}.

\subsection{Conjectured Combinatorial Formula}
\label{subsec:conj-rat-qcat}

We can specialize the rational $q,t$-Catalan numbers
described in~\cite[\S7]{LW-ptnid} to give a conjectured combinatorial 
interpretation for $\Cat_{\s,\rr}(q)$ for any $\s,\rr$ with $\gcd(\s,\rr)=1$.  
We also introduce a related, non-standard combinatorial interpretation
for the $q$-binomial coefficient $\tqbin{\s+\rr}{\s,\rr}{q}$ that
does not require the hypothesis $\gcd(\s,\rr)=1$.

To state these conjectures, fix $\s,\rr\in\NN$.  
Let $\ptnD(\s,\rr)$ be the
set of integer partitions whose diagrams fit in the triangle with
vertices $(0,0)$, $(0,\s)$, and $(\rr,\s)$.  More precisely, $\mu\in
\ptnD(\rr,\s)$ iff $\mu=(\mu_1,\mu_2,\ldots,\mu_\s)$ where
$\mu_1\geq\mu_2\geq\cdots\geq \mu_\s\geq 0$ and $\mu_j\leq \rr(\s-j)/\s$
for all $j$.
For example, the figure below shows that $\mu=(6,3,2,0,0)\in 
\ptnD(5,8)$.

\begin{center}
\epsfig{file=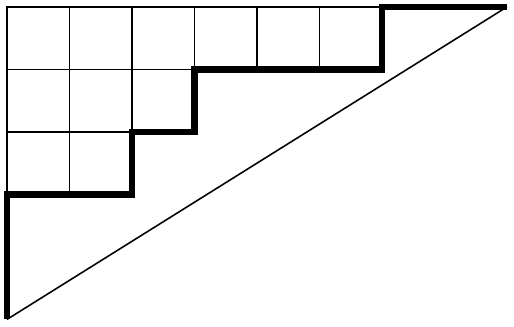,scale=0.7}
\end{center}

For each cell $c$ in the diagram of $\mu$, 
let $\arm(c)$ be the number of cells in the same row to the right of $c$,
and let $\leg(c)$ be the number of cells in the same column below $c$.
Let $h_{\rr,\s}^+(\mu)$ be the number of cells $c$ in the diagram of $\mu$
such that $-\s<\s\cdot\arm(c)-\rr\cdot\leg(c)\leq \rr$.
Let $h_{\rr,\s}^-(\mu)$ be the number of cells $c$ in the diagram of $\mu$
such that $-\s\leq \s\cdot\arm(c)-\rr\cdot\leg(c)<\rr$.
As above, $|\mu|$ is the total number of cells in the diagram of $\mu$.  
For the example partition shown above, 
we compute $|\mu|=11$ and $h_{8,5}^+(\mu)=h_{8,5}^-(\mu)=9$.
(In~\cite{LW-ptnid}, $h_{\rr,\s}^+$ and $h_{\rr,\s}^-$
 are denoted $h_{\rr/\s}^+$ and $h_{\rr/\s}^-$, respectively.) 

\begin{conj}\label{conj:rat-qcat} 
For all $\s,\rr\in\NN$ with $\gcd(\s,\rr)=1$, 
\[ \Cat_{\s,\rr}(q)=\sum_{\mu\in \ptnD(\s,\rr)} q^{|\mu|+h^+_{\rr,\s}(\mu)}. \]
\end{conj}

\begin{remark}
  We can use either $h^+_{\rr,\s}$ or $h^-_{\rr,\s}$ in this
  conjecture since the condition $\gcd(\s,\rr)=1$ ensures that the
  ``critical'' values $-\s$ and $\rr$ for
  $\s\cdot\arm(c)-\rr\cdot\leg(c)$ cannot actually occur.  The table
  below illustrates how each partition in $\ptnD(3,5)$ contributes to
  $\Cat_{3,5}(q,t) = 1 + q^2 + q^3 + q^4 + q^5 + q^6 + q^8$.
\end{remark}

\begin{table}[h]
  \centering
      \begin{tabular}{c>{$}c<{$}>{$}c<{$}>{$}c<{$}>{$}c<{$}}
$\mu$ & \mbox{frontier of }\mu & |\mu| & h^+_{5,3} & |\mu| + h^+_{5,3}
\\\toprule
$(3,1,0)$ &  NENEENEE & 4 & 4 & 8\\
$(2,1,0)$ &  NENENEEE & 3 & 3 & 6\\
$(3,0,0)$ &  NNEEENEE & 3 & 2 & 5\\
$(2,0,0)$ &  NNEENEEE & 2 & 2 & 4\\
$(1,1,0)$ &  NENNEEEE & 2 & 1 & 3\\
$(1,0,0)$ &  NNENEEEE & 1 & 1 & 2\\
$(0,0,0)$ &  NNNEEEEE & 0 & 0 & 0\\\bottomrule
  \end{tabular}
\end{table}

Next we describe our non-standard combinatorial interpretation for
$q$-binomial coefficients. Given $\mu\in \ptnR(\s,\rr)$, the
\textbf{frontier} of $\mu$ is the lattice path from $(0,0)$ to $(\rr,\s)$
obtained by taking north and east steps along the boundary of the
diagram of $\mu$.  For fixed $\s,\rr\in\NN$, assign the \textbf{level}
$\rr y-\s x$ to the lattice point $(x,y)\in\NN^2$. 
(These are the same levels for lattice points used 
in the proof of Proposition~\ref{prop:multinomial}.) 
For example, the partition shown above 
has frontier $NNEENENEEENEE$, which visits lattice points with levels
\[ 0,8,16,11,6,14,9,17,12,7,2,10,5,0. \]
Let $\ml_{\rr,\s}(\mu)$ be the minimum level appearing on the frontier of $\mu$.
Note that $\ml_{\rr,\s}(\mu)\leq 0$, with equality iff 
$\mu\in \ptnD(\s,\rr)$.

\begin{conj}\label{conj:nonstd-qbin} For all $\s,\rr\in\NN$,
\begin{equation}\label{eq:nonstd-qbin}
 \dqbin{\s+\rr}{\s,\rr}{q}
=\sum_{\mu\in \ptnR(\s,\rr)} q^{|\mu|+\ml_{\rr,\s}(\mu)+h^+_{\rr,\s}(\mu)}  
=\sum_{\mu\in \ptnR(\s,\rr)} q^{|\mu|+\ml_{\rr,\s}(\mu)+h^-_{\rr,\s}(\mu)}.
\end{equation}
\end{conj}

\begin{example}
  We compute the middle expression in Conjecture~\ref{conj:nonstd-qbin} for $a=2$, $b=3$.
  \begin{table}[h]
  \centering
  \begin{tabular}{c>{$}c<{$}>{$}c<{$}>{$}c<{$}>{$}c<{$}>{$}c<{$}}
$\mu$ & \mbox{frontier of }\mu & |\mu| & \ml_{3,2}(\mu) & h^+_{3,2} 
& |\mu| + \ml_{3,2}(\mu) + h^+_{3,2}\\\toprule
$(0,0)$  & NNEEE & 0 & 0 & 0 & 0\\
$(1,0)$  & NENEE & 1 & 0 & 1 & 2\\
$(2,0)$  & NEENE & 2 & -1 & 2& 3\\
$(3,0)$  & NEEEN & 3 & -3 & 2& 2\\
$(1,1)$  & ENNEE & 2 & -2 & 1& 1\\
$(2,1)$  & ENENE & 3 & -2 & 3& 4\\
$(3,1)$  & ENEEN & 4 & -3 & 4& 5\\
$(2,2)$  & EENNE & 4 & -4 & 3& 3\\
$(3,2)$  & EENEN & 5 & -4 & 5& 6\\
$(3,3)$  & EEENN & 6 & -6 & 4& 4\\\bottomrule
  \end{tabular}
  \end{table}
  The corresponding polynomial in $q$ is 
\[ \dqbin{2 + 3}{2,3}{q} = \frac{[5]_q[4]_q}{[2]_q} 
= \frac{1-q^5}{1-q}\cdot\frac{1-q^4}{1-q}\cdot\frac{1-q}{1-q^2} 
= 1+q+2q^2+2q^3+2q^4+q^5+q^6. \]
\end{example}

We have confirmed Conjectures~\ref{conj:rat-qcat}
and~\ref{conj:nonstd-qbin} for $a,b\leq 12$ using a Sage~\cite{sage}
worksheet~\cite{sage-worksheet}.  The rest of this section is devoted
to a bijective proof of the following theorem.

\begin{theorem}\label{thm:ratcat}
Conjecture~\ref{conj:rat-qcat} holds for given coprime $a,b\in\NN$
if and only if Conjecture~\ref{conj:nonstd-qbin} holds for this choice
of $a,b$.
\end{theorem}

\subsection{Alternate Description of $h^+$ and $h^-$}
\label{subsec:alt-stat}

We begin with a lemma showing how to compute $h^+_{\rr,\s}(\mu)$
and $h^-_{\rr,\s}(\mu)$ from the levels on the frontier of $\mu$.

\begin{lemma}\label{lem:h-via-labels}
Given $\s,\rr\in\NN$ and $\mu\in \ptnR(\s,\rr)$, 
let the frontier of $\mu$ be $w_1w_2\cdots w_{\s+\rr}\in\{N,E\}^{\s+\rr}$,
and let the levels on the frontier
of $\mu$ be $l_0=0,l_1,l_2,\ldots,l_{\s+\rr}=0$.
\begin{itemize}
\item[(a)] $h^+_{\rr,\s}(\mu)$ is the number of pairs $i<j$
with $w_i=E$, $w_j=N$ and $1\leq l_{i-1}-l_{j-1}\leq \s+\rr$.
\item[(b)] $h^-_{\rr,\s}(\mu)$ is the number of pairs $i<j$ 
with $w_i=E$, $w_j=N$, and $1\leq l_j-l_i\leq \s+\rr$.  
\end{itemize}
\end{lemma}
\begin{proof}
(a) First note that there is a bijection between the set of cells $c$
in the diagram of $\mu$ and the set of pairs $i<j$ with $w_i=E$ and $w_j=N$.
This bijection maps $c$ to the unique pair $(i,j)$ such that 
$w_i$ is the east step on the frontier due south of $c$, and
$w_j$ is the north step on the frontier due east of $c$.
We can use the coordinates of the starting vertices of these steps
to determine whether $c$ contributes to $h^+_{\rr,\s}(\mu)$:
Suppose $w_i$ goes from $(u,v)$ to $(u+1,v)$. By definition of
$\arm(c)$ and $\leg(c)$, we see that $w_j$ must go from
$(u+\arm(c)+1,v+\leg(c))$ to $(u+\arm(c)+1,v+\leg(c)+1)$. 
Using these coordinates, we have $l_{i-1}=\rr v-\s u$ and
$l_{j-1}=\rr(v+\leg(c))-\s(u+\arm(c)+1)$, so
$l_{i-1}-l_{j-1}=\s+\s\cdot\arm(c)-\rr\cdot\leg(c)$. 
This cell $c$ contributes to $h^+_{\rr,\s}(\mu)$ iff
$-\s<\s\cdot\arm(c)-\rr\cdot\leg(c)\leq \rr$ iff $0<l_{i-1}-l_{j-1}\leq \s+\rr$.

(b) Keep the notation from (a). We find that
$l_i=\rr v-\s(u+1)$ and $l_j=\rr(v+\leg(c)+1)-\s(u+\arm(c)+1)$, so
$l_j-l_i=\rr\cdot\leg(c)-\s\cdot\arm(c)+\rr$.
The cell $c$ contributes to $h^-_{\rr,\s}(\mu)$ iff
$-\s\leq \s\cdot\arm(c)-\rr\cdot\leg(c)<\rr$ iff
$-\s-\rr\leq \s\cdot\arm(c)-\rr\cdot\leg(c)-\rr<0$ iff $0<l_j-l_i\leq \s+\rr$.  
\end{proof}

Henceforth, we discuss only the statistic $h^+_{\rr,\s}$;
the arguments needed for $h^-_{\rr,\s}$ are analogous.

\subsection{Cyclic Shift Analysis}
\label{subsec:cyclic-shift}

Define a \emph{cyclic-shift map} $C$ on words by letting
$C(w_1w_2\cdots w_{\s+\rr})=w_2\cdots w_{\s+\rr}w_1$. We get an
associated bijection $C$ on $\ptnR(\s,\rr)$ as follows. 
Given $\mu\in \ptnR(\s,\rr)$, let $w$ be the frontier of $\mu$; then
let $C(\mu)$ be the unique partition in $\ptnR(\s,\rr)$ with frontier $C(w)$.
The next lemma reveals how the cyclic-shift map affects $h^+_{\rr,\s}$.

\begin{lemma}\label{lem:cyc-shift}
Given $\s,\rr\in\NN$ and $\mu\in \ptnR(\s,\rr)$, 
let the frontier of $\mu$ be $w_1w_2\cdots w_{\s+\rr}\in\{N,E\}^{\s+\rr}$,
with corresponding vertex levels $l_0,\ldots,l_{\s+\rr}$.  
Let $\Delta h^+=h^+_{\rr,\s}(C(\mu))-h^+_{\rr,\s}(\mu)$.
\begin{itemize}
\item[(a)] If $w_1=N$, then 
\[ \Delta h^+= |\{i>0: w_i=E\mbox{ and }1\leq l_{i-1}\leq \s+\rr\}|
             = |\{k>0: 1\leq l_{k-1}\leq \rr\}|. \]
\item[(b)] If $w_1=E$, then 
\[ \Delta h^+=-|\{j>0: w_j=N\mbox{ and }1\leq -l_{j-1}\leq \s+\rr\}|
             =-|\{k>0: 1\leq -l_{k-1}\leq \s\}|. \] 
\end{itemize} 
\end{lemma}
\begin{proof}
(a) Assume $w_1=N$. The frontier of $C(\mu)$ is $w_1'w_2'\cdots w_{\s+\rr}'$,
where $w_{\s+\rr}'=N$ and $w_i'=w_{i+1}$ for $1\leq i<\s+\rr$. The levels of
the vertices on the frontier of $C(\mu)$ are $l_0',l_1',\ldots,l_{\s+\rr}',$
where $l_{\s+\rr}'=0$ and $l_i'=l_{i+1}-\rr$ for $0\leq i<\s+\rr$.
Let us use Lemma~\ref{lem:h-via-labels} to compute
$h^+_{\rr,\s}(\mu)$ and $h^+_{\rr,\s}(C(\mu))$. Note that $(i,j)$ satisfies
$0<i<j<\s+\rr$, $w_i'=E$, $w_j'=N$, and $1\leq l_{i-1}'-l_{j-1}'\leq \s+\rr$
iff $(i+1,j+1)$ satisfies $1<i+1<j+1\leq \s+\rr$, $w_{i+1}=E$, $w_{j+1}=N$, and
$1\leq l_i-l_j\leq \s+\rr$. So the contributions to $h^+_{\rr,\s}$
from these pairs will cancel in the computation of $\Delta h^+$.
We still need to consider pairs $(i,j)$ contributing to $h^+_{\rr,\s}(C(\mu))$
where $j=\s+\rr$, and pairs $(i,j)$ contributing to $h^+_{\rr,\s}(\mu)$
where $i=1$. But, since $w_1=N$, no pair of the second type causes a contribution to
$h^+_{\rr,\s}(\mu)$. 
Since $w_{\s+\rr}'=N$ and $l_{\s+\rr-1}'=-\rr$,
a pair $(i,\s+\rr)$ contributes to $h^+_{\rr,\s}(C(\mu))$ iff
$w_i'=E$ and $1\leq l_{i-1}'+\rr\leq \s+\rr$ iff 
$w_{i+1}=E$ and $1\leq l_i\leq \s+\rr$. Replacing $i$ by $i-1$ 
gives the first formula in (a).

To prove the second formula in (a), let $A$ be the set of lattice points
on the frontier of $\mu$ with levels in $[1,\s+\rr]$
and followed by an east step, and let $B$ be the set of lattice points
on the frontier of $\mu$ with levels in $[1,\rr]$.
It suffices to define a bijection $f:A\rightarrow B$.
Given $(x,y)\in A$ with level $i\in [1,\s+\rr]$, 
the next point on the frontier of $\mu$ is $(x+1,y)$,
which has level $i-\s\in [1-\s,\rr]$. If this level is positive,
then set $f(x,y)=(x+1,y)\in B$. Otherwise, keep stepping
forward along the frontier (wrapping back at the end from 
$(\rr,\s)$ to $(0,0)$ if needed)
until a positive level is reached, and let $f(x,y)$ be
the location of this level. Since north steps increase
the level value in increments of $\rr$, the positive 
level reached must lie in $[1,\rr]$. Also, such a level
must exist, as $(x,y)$ has a positive level. The inverse 
bijection is similar: given $(u,v)\in B$, scan backwards
along the frontier until one hits a point with level in $[1,\s+\rr]$
that is reached by following an east step backwards.
One sees that $f(x,y)$ is always the first point after $(x,y)\in A$
with positive level, whereas $g(u,v)$ is always the first point
before $(u,v)\in B$ with positive level. It follows that
$f\circ g$ and $g\circ f$ are identity maps.  

(b) Keep the notation from (a). In this case, since $w_1=E$, the new 
levels $l_0',\ldots,l_{\s+\rr}'$ satisfy $l_{\s+\rr}'=0$
and $l_i'=l_{i+1}+\s$ for $0\leq i<\s+\rr$.  
As in (a), a pair $(i,j)$ with $w_i'=E$ and $w_j'=N$ and $0<i<j<\s+\rr$ 
contributes to $h^+_{\rr,\s}(C(\mu))$ iff the pair $(i+1,j+1)$ contributes 
to $h^+_{\rr,\s}(\mu)$. Since $w_{\s+\rr}'=w_1=E$, no pairs
$(i,j)$ with $j=\s+\rr$ contribute to $h^+_{\rr,\s}(C(\mu))$.
On the other hand, a pair $(1,j)$ contributes to $h^+_{\rr,\s}(\mu)$
iff $w_j=N$ and $1\leq l_0-l_{j-1}\leq \s+\rr$.
So $-\Delta h^+$ is the number of $j>0$ with $w_j=N$ and 
$1\leq -l_{j-1}\leq \s+\rr$.

To prove the second formula in (b), let $A$ be the set of lattice
points on the frontier of $\mu$ with levels in $[-(\s+\rr),-1]$ and
followed by a north step, and let $B$ be the set of lattice points on
the frontier of $\mu$ with levels in $[-\s,-1]$. Define a bijection
$f:A\rightarrow B$ by letting $f(x,y)$ be the first point following
$(x,y)\in A$ along the frontier that has a negative level. The inverse
bijection $g:B\rightarrow A$ sends $(u,v)\in B$ to the first point
preceding $(u,v)$ along the frontier that has a negative level.  Both
maps ``wrap around'' at $(\rr,\s)$ and $(0,0)$ when needed.  One checks
that $f$ maps $A$ into $B$, $g$ maps $B$ into $A$, and the two maps
are inverses.
\end{proof}

\subsection{Cyclic Shift Orbits when $\gcd(\s,\rr)=1$}
\label{subsec:cyc-shift-orbits}

We now impose the hypothesis $\gcd(\s,\rr)=1$.
In this case, the argument in Step 2 of the proof of 
Proposition~\ref{prop:multinomial} shows that the levels $l_0,l_1,\ldots,
l_{a+b-1}$ on any lattice path from $(0,0)$ to $(\rr,\s)$ are all 
\emph{distinct}.  
Because each shifted path has a different minimum level, it follows that the $\s+\rr$ possible cyclic shifts of this lattice path are 
all distinct.
Moreover, exactly one of the $\s+\rr$ cyclic shifts has minimum level zero,
and all other cyclic shifts have a negative minimum level
(see the proof of Proposition~\ref{prop:multinomial} and \cite[\S12.1]{loehr-bij}).

We can rephrase these results in terms of partitions.  
Given $\mu\in \ptnR(\s,\rr)$, there exists a unique partition 
$\mu^0\in \ptnD(\s,\rr)$ that can be obtained by cyclically shifting the 
frontier of $\mu$.  The $\s+\rr$ partitions $\mu^0,C(\mu^0),\ldots,
C^{\s+\rr-1}(\mu^0)$ are all distinct.
By induction on the number of cyclic shifts needed to pass from
$\mu^0$ to $\mu$, one may check that $|\mu^0|=|\mu|+\ml_{\rr,\s}(\mu)$. 
The following lemma is the final ingredient needed for the proof
of Theorem~\ref{thm:ratcat}.

\begin{lemma}\label{lem3} 
Assume $\gcd(\s,\rr)=1$.
Given $\mu^0\in \ptnD(\s,\rr)$, let $\mu^0,\mu^1,\ldots,\mu^{\s+\rr-1}$
be the $\s+\rr$ distinct partitions obtained by cyclically
shifting the frontier of $\mu^0$. Then
\[ \sum_{i=0}^{\s+\rr-1} q^{h^+_{\rr,\s}(\mu^i)}
  =q^{h^+_{\rr,\s}(\mu^0)}[\s+\rr]_q. \]
More precisely, suppose the $\s+\rr$ (distinct) levels on 
the frontier of $\mu^0$ are $0=i_0<i_1<i_2<\cdots<i_{\s+\rr-1}$.
For $0\leq k<\s+\rr$, if $\mu^k$ is the cyclic shift of $\mu^0$ with 
minimum level $-i_k$, then $h^+_{\rr,\s}(\mu^k)-h^+_{\rr,\s}(\mu^0)=k$.
\end{lemma}
\begin{proof}
Fix a partition $\mu$ in the orbit of $\mu^0$ with minimum level $-i_k$.
Let $\nu=C(\mu)$, and consider two cases.

Case~1: The frontier of $\mu$ begins with a north step.
Then $\nu$ has minimum level $-i_k-\rr$. 
In the diagram for $\mu$,
the lattice point $(0,0)$ has level zero. This lattice point was
moved to the origin by cyclically shifting some lattice point
$(x,y)$ on the frontier of $\mu^0$. Since $\ml_{\rr,\s}(\mu)=-i_k$
and $\ml_{\rr,\s}(\mu^0)=0$, it follows that $(x,y)$ has level $i_k$.
Similarly, the point $(0,0)$ in the diagram of $\nu$ arose by cyclically
shifting some point $(u,v)$ on the frontier of $\mu^0$ with level $i_k+\rr$,
where $i_k+\rr=i_j$ for some $j>k$.  Now, keeping in mind that all levels
are distinct, $j-k$ is precisely the number of 
levels on the frontier of $\mu^0$ whose values are between $i_k$ and $i_k+\rr$,
excluding $i_k$ and including $i_k+\rr$. Cyclically shifting $\mu^0$ to $\mu$,
this means that $j-k$ is the number of levels on the frontier of $\mu$
with values in the interval $[1,\rr]$. By Lemma~\ref{lem:cyc-shift}(a), we have
$h^+_{\rr,\s}(\nu)-h^+_{\rr,\s}(\mu)=j-k$.

Case~2:  The frontier of $\mu$ begins with an east step.
Then $\nu$ has minimum level $-i_k+\s$. As in Case 1,
the origin in the diagram of $\mu$ is the cyclic shift of a point
$(x,y)$ on the frontier of $\mu^0$ at level $i_k$, whereas
the origin in the diagram of $\nu$ is the cyclic shift of a point
$(u,v)$ on the frontier of $\mu^0$ at level $i_k-\s$.
We know $i_k-\s=i_j$ for some $j<k$. Here, $|j-k|$ is
the number of levels in $\mu^0$ with values between $i_j=i_k-\s$ (inclusive)
and $i_k$ (exclusive). Cycling $\mu^0$ to $\mu$, we see that
$|j-k|$ is the number of levels in $\mu$ with values in the interval $[-\s,-1]$.
By Lemma~\ref{lem:cyc-shift}(b),
$h^+_{\rr,\s}(\nu)-h^+_{\rr,\s}(\mu)=-|j-k|=j-k$.

Repeatedly using the two cases, the claimed formula
$h^+_{\rr,\s}(\mu)-h^+_{\rr,\s}(\mu^0)=k$ now follows 
by induction on the number of cyclic shift steps required to
go from $\mu^0$ to $\mu$.
\end{proof}

We know $\ptnR(\s,\rr)$ is the disjoint union of the orbits of the cyclic-shift
map $C$, where each orbit has a unique representative 
$\mu^0\in \ptnD(\s,\rr)$.
Combining this fact with the last lemma, we see that
\[ \sum_{\mu\in \ptnR(\s,\rr)}q^{|\mu|+\ml_{\rr,\s}(\mu)+h_{\rr,\s}^+(\mu)}
 =[\s+\rr]_q\sum_{\mu^0\in \ptnD(\s,\rr)}q^{|\mu^0|+h_{\rr,\s}^+(\mu^0)}. \]
This equation shows the equivalence of Conjecture~\ref{conj:rat-qcat}
and Conjecture~\ref{conj:nonstd-qbin} for each fixed $\s,\rr$ 
with $\gcd(\s,\rr)=1$. 
So the proof of Theorem~\ref{thm:ratcat} is complete.

Below we illustrate Lemma~\ref{lem3} for $\mu^0 = (2,1,0) \in D(3,5)$.
Cells contributing to $h^{+}_{5,3}$ are labeled with the corresponding
value of $l_{i-1}-l_{j-1}$ as explained in
Lemma~\ref{lem:h-via-labels}.  The vertex labels contributing to
$\Delta h^{+}$ are italicized and written in blue (red) for paths
starting with an east (north) step.
\begin{center}
\includegraphics[scale=.4]{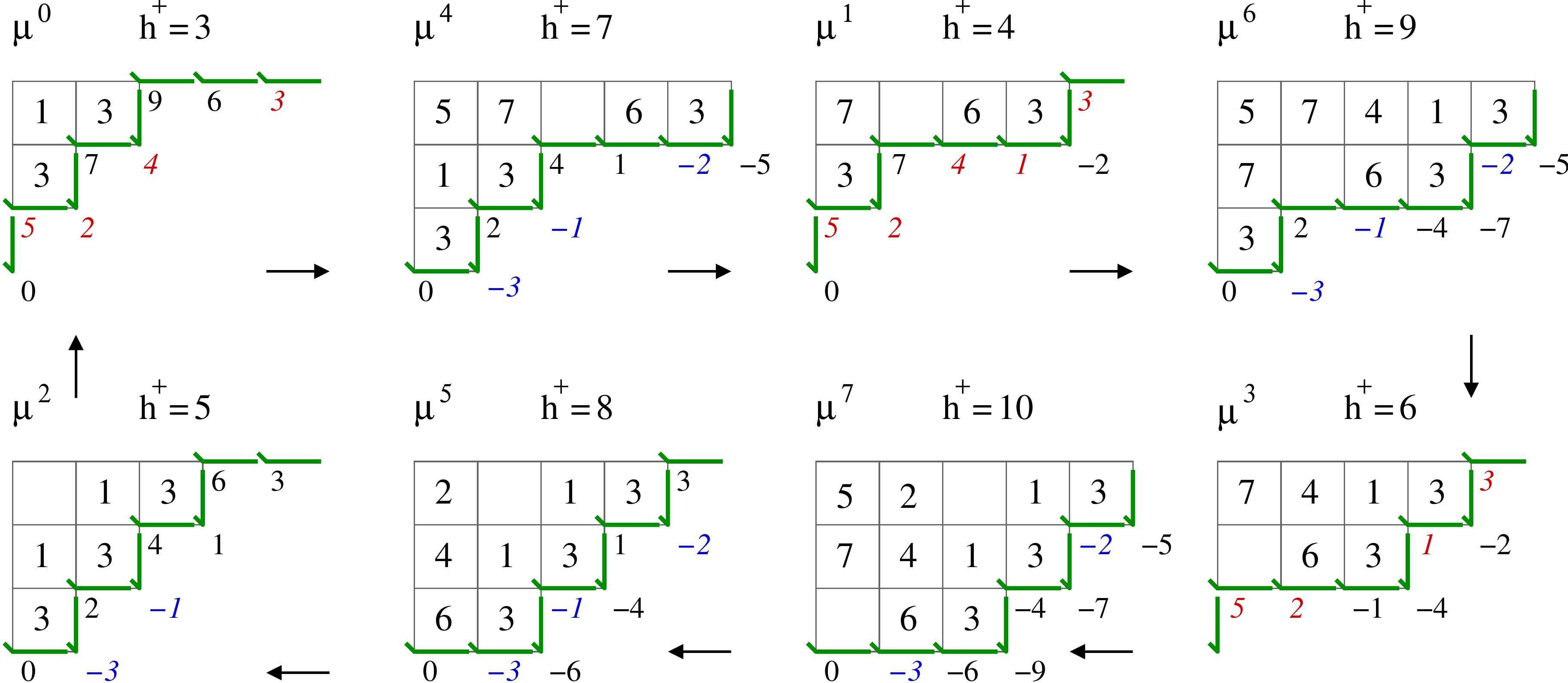}
\end{center}

\subsection{Remarks on Conjecture~\ref{conj:nonstd-qbin}.}
\label{subsec:conclusion}

A natural approach to proving Conjecture~\ref{conj:nonstd-qbin}
is to show that the combinatorial formulas in~\eqref{eq:nonstd-qbin}
satisfy the recursions~\eqref{eq:qbin-rec}. One difficulty here
is determining how to divide $\ptnR(\s,\rr)$ into two disjoint sets that
are $q$-counted by the two terms in these recursions. The approach
used for the standard interpretation of $q$-binomial coefficients
(removing the first step or the last step of the frontier) does not work.  
Another major obstacle is the fact that when the dimension of the rectangle 
is reduced from $(\rr,\s)$ to $(\rr,\s-1)$ or $(\rr-1,\s)$, the statistic
changes from $h_{\rr,\s}^+$ to $h_{\rr,\s-1}^+$ or $h_{\rr-1,\s}^+$, so that
the power of $q$ for a particular partition changes unpredictably.  

\section{The Shuffle Conjecture}
\label{sec:shuffle}

To prepare for our discussion of rational $q,t$-parking functions
in the next section, this section reviews some combinatorial conjectures
for the Hilbert series and Frobenius series of $\DR_n$ involving
classical $q,t$-parking functions. These conjectures first appeared
in~\cite{HHLRU,HL-park}.

\subsection{The $\dinv$ Statistic.}
\label{subsec:dinv}

Let $P\in\PF_n$ be a labeled Dyck path.
We have already defined $\area(P)$ in \S\ref{subsec:graded-pfn}.
We now define a second statistic $\dinv(P)$ (the \textbf{diagonal
inversion count} of $P$) as follows. Use $P$ to
define vectors $(g_1,\ldots,g_n)$ and $(p_1,\ldots,p_n)$
where $g_i$ is the number of area cells in the $i$-th row
from the bottom, and $p_i$ is the label in the $i$-th row
from the bottom. For example, the labeled Dyck path shown
at the end of \S\ref{subsec:PF} has $(g_1,\ldots,g_5)=(0,1,1,2,1)$
and $(p_1,\ldots,p_5)=(3,5,1,4,2)$. Let $\dinv(P)$
be the number of pairs $i<j$ such that either $g_i=g_j$ and $p_i<p_j$,
or $g_i=g_j+1$ and $p_i>p_j$. Our example object has $\dinv(P)=2$,
corresponding to the pairs $(i,j)=(3,5)$ and $(i,j)=(4,5)$.
Haglund, Haiman, and Loehr conjectured the following formula
for the Hilbert series of the doubly-graded diagonal coinvariant module.

\begin{conj}[\cite{HL-park}]
For all $n>0$, 
$$\Hilb_{\DR_n}(q,t)=\sum_{P\in\PF_n} q^{\area(P)}t^{\dinv(P)}.$$
\end{conj}

\subsection{The $\zetamap$ Map.}
\label{subsec:zeta}

The unusual statistic $\dinv$ is actually a disguised version of
the natural statistic $\area$. To explain this, we recall the definition
of a bijection $\zetamap$ between two different
encodings of parking functions, which are two different ways to label
Dyck paths. We are already familiar with the first labeling scheme
(filling the vertical runs with increasing labels). Call this the 
{\bf coset notation}. 
A parking function under the second labeling scheme will also consist
of a Dyck path $D$ paired with a permutation.  In this case, the
entries of the permutation will be written along the diagonal $y=x$.
However, just as the labels in the coset notation must increase along
each vertical run, there are restrictions on the permutation here as
well.  To describe these restrictions, we label each square above
$y=x$ by the ordered pair $(i,j)$ where $i$ is the permutation entry
in the same column and $j$ is the permutation entry in the same row.
The pair consisting of the Dyck path and the permutation 
is a parking function iff
each square lying both immediately above an east step of $D$ and
immediately to the left of a north step of $D$ has a label $(i,j)$
satisfying $i < j$. (We call such squares \textbf{left-turns} of $D$.)
The figure below shows an order-$5$
parking function using this new labeling system.
Call this the {\bf root notation} for parking functions. 
\bigskip
\begin{center}
  \includegraphics[scale=.6]{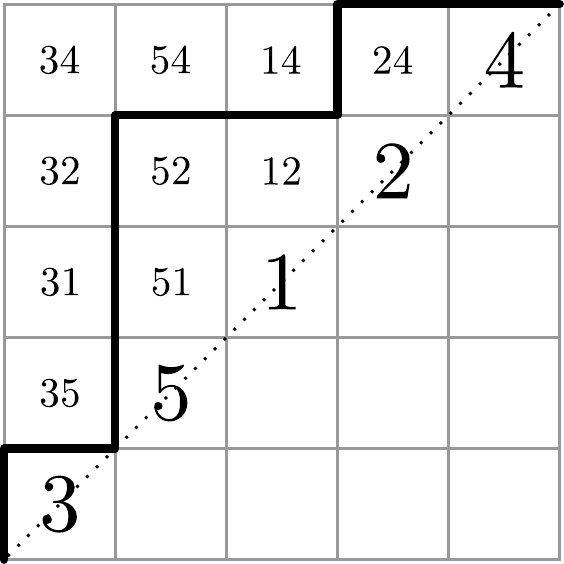}\label{fig:otherway}
\end{center} 
Note that the pairs labeling the left-turns are $35$ and $14$, which are
indeed increasing.  We remark that parking functions in root notation
are in natural bijection with regions of the Shi hyperplane
arrangement (see Armstrong~\cite{armstrong}).  

To define the map $\zetamap$, begin with a parking function
$P\in\PF_n$ in coset notation. Define the \textbf{diagonal reading
word} of $P$, denoted $\drw(P)$, by reading labels along diagonals of
slope $1$, working from higher diagonals to lower diagonals,
and scanning each diagonal from northeast to southwest.
For example, the $P$ displayed in~\S\ref{subsec:graded-pfn} 
(in coset notation) has $\drw(P)=42153$.
Let $l_1l_2\cdots l_n$ be the reverse of $\drw(P)$.  
Define a set $\mathsf{valleys}(P)$ consisting of the
ordered pairs $(j,k)$ such that the label $l_j$ occurs in the box
just below $l_k$ in coset notation. Finally, let $\zetamap(P)$ be
the parking function in root notation with 
\begin{itemize}
\item word $l_1l_2\cdots l_n$ on the diagonal and
\item Dyck path having
left turns at the squares containing the labels
$\mathsf{valleys}(P)$.
\end{itemize}
One may check that the elements of $\mathsf{valleys}(P)$ are
``non-nesting'' so that such a path exists. The following figure gives
an example of this map: \bigskip
\begin{center}
\includegraphics[scale=.6]{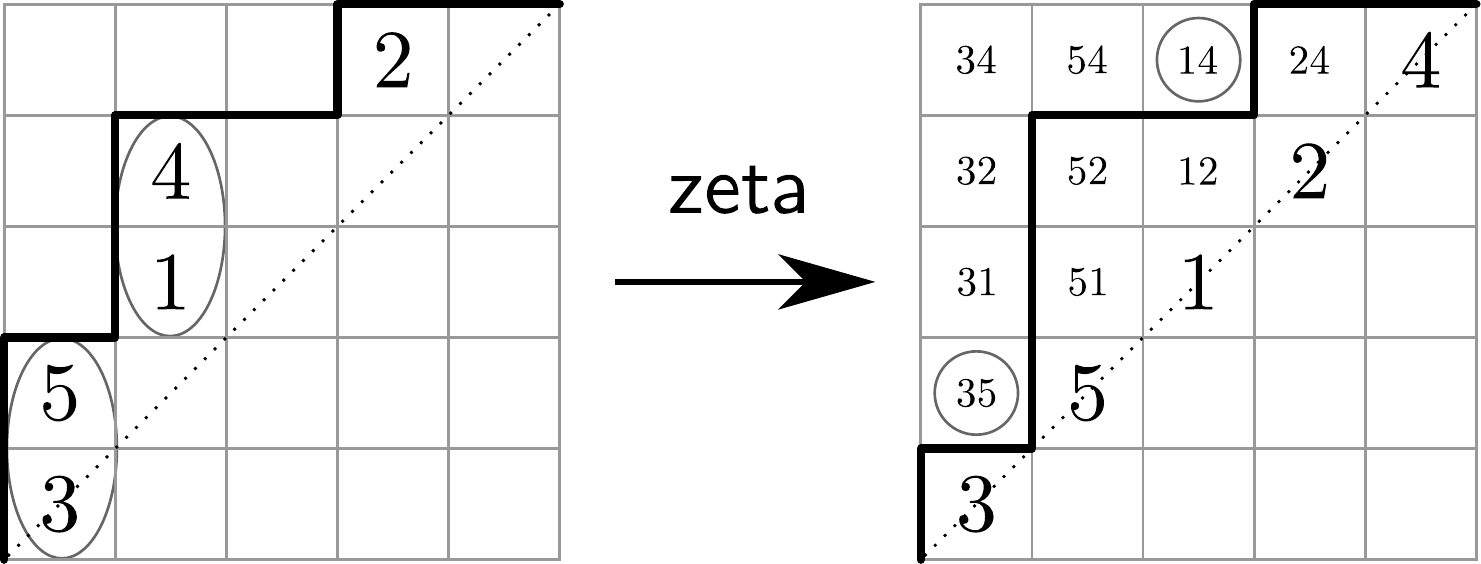}
\end{center}

\begin{remark}
  We regard the coset notation and the root notation as analogous to
  the two combinatorial ways to view an affine Weyl group --- via
  minimal coset representatives for the finite Weyl group, and via the
  semi-direct product of the finite Weyl group and the root
  lattice. For more details, see~\cite{armstrong}.
\end{remark}

\subsection{The $\areap$ Statistic.}
\label{subsec:areap}

Next, we define the $\areap$ statistic. If $P\in\PF_n$ is a parking
function in coset notation with corresponding parking function
$\zetamap(P)$ in root notation, we let $\areap(\zetamap(P))$ be the
number of boxes fully between the path and the diagonal in
$\zetamap(P)$ that are labeled by {\em increasing pairs}. For our
running example $P=(2,4,1,2,1)$, we have $\areap(\zetamap(P))=2$,
coming from the boxes labeled $12$ and $24$.  It can be shown
that $\dinv(P)=\areap(\zetamap(P))$ for all $P\in\PF_n$
(cf.~\cite{HL-park}, which proves this using a different description
of the map $\zetamap$). Accordingly, we can restate the Hilbert
series conjecture as follows.

\begin{conj}[\cite{HL-park}]
For all $n>0$,
\[ \Hilb_{\DR_n}(q,t)=\sum_{P\in\PF_n} q^{\area(P)}t^{\areap(\zetamap(P))}.\]
\end{conj}

\subsection{The Shuffle Conjecture}
\label{subsec:shuffle}

We need two more ingredients to state the Shuffle
Conjecture~\cite{HHLRU} giving a combinatorial formula for the
Frobenius series of the doubly-graded module $\DR_n$.  First, we say
that a formal power series in $\CC[[\x]]$ is {\bf quasisymmetric} iff
it is invariant under shifts of the variables that preserve the order
of indices.  Let $\QSym\subseteq \CC[[\x]]$ denote the subalgebra of
quasisymmetric functions. Gessel's {\bf fundamental basis} of $\QSym$
is indexed by pairs $n,S$ where $S\subseteq [n-1]=\{1,2,\ldots,n-1\}$.
This basis is defined by
\begin{equation*}
F_{n,S}(\x)= \sum_{\substack{1\leq i_1\leq i_2\leq \cdots \leq i_n
\\ i_j<i_{j+1} \text{ if } j\in S}} x_{i_1}x_{i_2}\cdots x_{i_n}.
\end{equation*} 
Second, given a parking function $P\in\PF_n$ with diagonal reading word
$\drw(P)$, define the \textbf{inverse descent set}
$\IDes(P)\subseteq [n-1]$ by letting $j\in\IDes(P)$ if and only if
$j+1$ occurs to the left of $j$ in $\drw(P)$, for all $j\in [n-1]$.  
Our running example has $\drw(P)=42153$ and $\IDes(P)=\{1,3\}$.
\begin{conj}[Shuffle Conjecture~\cite{HHLRU}]
\label{conj:shuffle}
  The bi-graded Frobenius characteristic of the ring of diagonal coinvariants 
  $\DR_n$ (considered as an $S_n$-module bi-graded by $\x$-degree and
  $\y$-degree) satisfies
\begin{equation}\label{eq:shuffle}
\Frob_{\DR_n}(\x;q,t)= 
\sum_{P\in\PF_n} q^{\area(P)} t^{\areap(\zetamap(P))} F_{n,\IDes(P)}(\x)=
\sum_{P\in\PF_n} q^{\area(P)} t^{\dinv(P)} F_{n,\IDes(P)}(\x).
\end{equation}
\end{conj}

Recall that \emph{increasing} parking functions (in coset notation)
correspond bijectively with \emph{unlabeled} Dyck paths. One
can check that when $D$ is increasing, $\dinv(D)$ is the number
of $i<j$ with $g_i-g_j\in\{0,1\}$, and
$\areap(\zetamap(D))=\area(\zetamap(D))$.
As a special case of the Shuffle Conjecture, the coefficient of the
sign representation in $\Frob_{\DR_n}$ is given by the 
\textbf{$q,t$-Catalan number}
\begin{equation}
\label{eq:qtcat}
  \Cat_n(q,t)=\sum_{D\in\patD(N^nE^n)} q^{\area(D)} t^{\area(\zetamap(D))}
             =\sum_{D\in\patD(N^nE^n)} q^{\area(D)} t^{\dinv(D)}.
\end{equation}
This special case of the conjecture has been proved by Garsia and
Haglund~\cite{qtcat-proof}.

The full Shuffle Conjecture is quite remarkable since it contains many
famous combinatorial polynomials at various specializations of $q$,
$t$, and $\x$. It has also been exceptionally difficult to prove. One
issue is that the combinatorial formula is experimentally seen to be 
symmetric in $q$ and $t$ (as it must be to satisfy the conjecture); however, 
this symmetry is not visible in the combinatorics. In fact, it is an open 
problem to prove combinatorially that the expression on the right side 
of~\eqref{eq:shuffle} is symmetric in $q$ and $t$.


\section{Rational $q,t$-Catalan Numbers}
\label{sec:qt-ratcat}

Throughout this section, fix coprime integers $a,b\in\NN$.
We now generalize the combinatorics of $q,t$-parking functions
and $q,t$-Catalan numbers from the classical parking functions
$\PF_{n,n+1}$ to the rational parking functions $\PF_{a,b}$
In doing so, we will define a rational
version of the shuffle conjecture~\eqref{eq:shuffle}.
Our new formula should be
the Frobenius series for some natural representation of $S_a$ with 
a symmetric bi-grading, generalizing the ring of diagonal
coinvariants. A conjecturally suitable representation of $S_a$ has
been constructed using the theory of rational Cherednik
algebras~\cite{OY}, but at present the construction is very
difficult. It is an open problem to construct this representation 
more directly.

\subsection{The Sweep Map}
\label{sec:sweep} 

Recall from~\eqref{eq:qtcat}
that the classical $q,t$-Catalan numbers of Garsia and
Haiman~\cite{GH-qtcat} can be defined in terms of classical Dyck paths
$\patD(N^n E^n)$ using one natural statistic (called $\area$) and one
strange map (called $\zetamap$). 
We already generalized the $\area$ statistic to rational Dyck paths
in \S\ref{subsec:PFab-graded}.  Next we define a map, called the
\textbf{sweep map}, that will play the role of $\zetamap$ in the rational
case. This map is studied in detail in~\cite{sweep}; see also~\cite{AHJ}.
  
\begin{definition}[The Sweep Map]
Consider a lattice path $w=w_1w_2\cdots w_{a+b}\in \patD(N^a E^b)$.
As in Step 2 of the proof of Proposition~\ref{prop:multinomial},
define \textbf{levels} $l_0,l_1,\ldots,l_{a+b}$ of lattice points on
the path by setting $l_0=0$ and, for $1\leq i\leq a+b$, 
setting $l_{i}=l_{i-1}+b$ if $w_{i-1}=N$ and $l_i =
l_{i-1}-a$ if $w_{i-1}=E$. For each step $w_i$ in the path,
let $l_{i-1}$ be the \textbf{wand label} associated to this step.
Define $\sweep(w)$ to be the path obtained by sorting the steps $w_i$ 
into increasing order according to their wand labels, and then 
erasing all the labels.
\end{definition}

For example, 
the figure below computes $\sweep(w)$ for the $(7,10)$-Dyck path
\begin{equation*}
  w = NENENNEENNNEEEEEE \in \patD(N^7 E^{10}).
\end{equation*}
Next to each step in $\sweep(w)$, we have written the wand label
of the corresponding step in $w$, although these labels are not
retained as part of the final output.
\begin{center}
      {\scalebox{.4}{\includegraphics{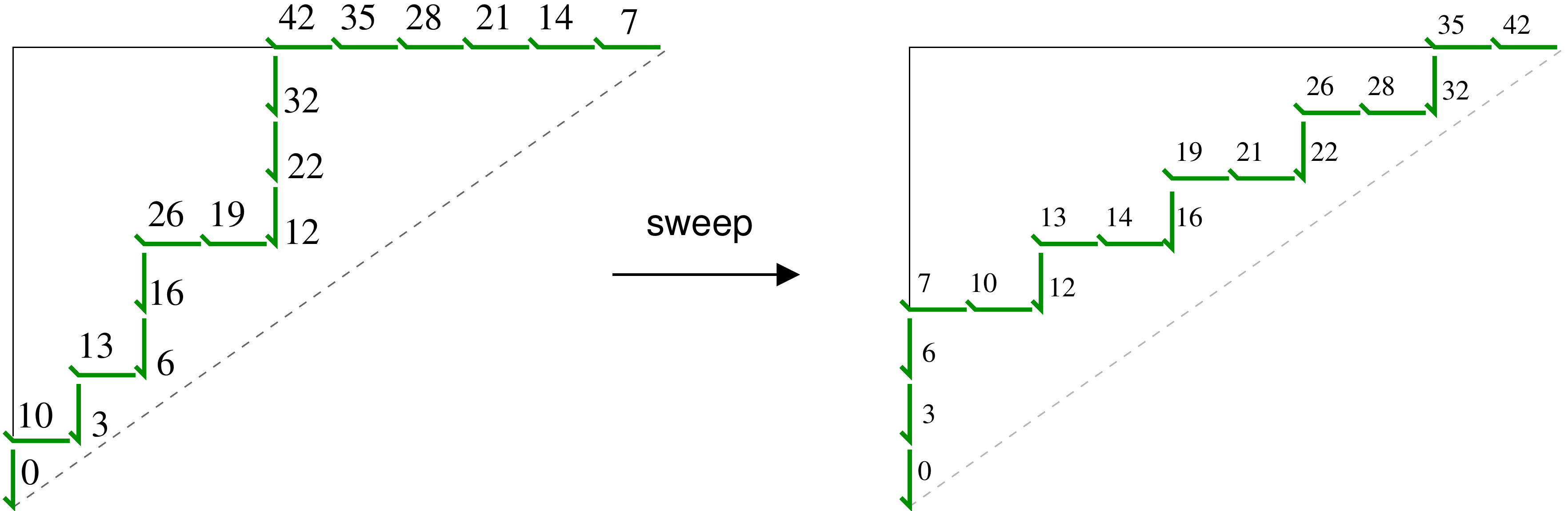}}}
\end{center}

Since $\sweep$ acts by reordering path steps, it is immediate from the
definition that if $D\in\patD(N^a E^b)$, then $\sweep(D)\in \patR(N^a
E^b)$. It is true but non-obvious that $\sweep(D)$ is actually an $(a,b)$-Dyck 
path; see~\cite{sweep} for a proof.  

\subsection{Rational $q,t$-Catalan Numbers} 
\label{subsec:ratqt}

We use the $\area$ statistic and
the $\sweep$ map to define the {\bf rational $q,t$-Catalan number}
\begin{equation}\label{eq:rat-qtcat}
  \Cat_{a,b}(q,t)=\sum_{D\in\patD(N^a E^b)} q^{\area(D)} t^{\area(\sweep(D))}.
\end{equation}
For example, we have
\begin{equation*}
\Cat_{5,8}(q,t)= \begin{bmatrix} 
. & . & . & . & . & . & . & . & . & . & . & . & . & . & 1 \\
. & . & . & . & . & . & . & . & . & . & 1 & 1 & 1 & 1 & . \\
. & . & . & . & . & . & . & . & 2 & 2 & 2 & 1 & 1 & . & . \\
. & . & . & . & . & . & 1 & 3 & 3 & 2 & 1 & 1 & . & . & . \\
. & . & . & . & . & 2 & 4 & 3 & 2 & 1 & 1 & . & . & . & . \\
. & . & . & . & 2 & 4 & 3 & 2 & 1 & 1 & . & . & . & . & . \\
. & . & . & 1 & 4 & 3 & 2 & 1 & 1 & . & . & . & . & . & . \\
. & . & . & 3 & 3 & 2 & 1 & 1 & . & . & . & . & . & . & . \\
. & . & 2 & 3 & 2 & 1 & 1 & . & . & . & . & . & . & . & . \\
. & . & 2 & 2 & 1 & 1 & . & . & . & . & . & . & . & . & . \\
. & 1 & 2 & 1 & 1 & . & . & . & . & . & . & . & . & . & . \\
. & 1 & 1 & 1 & . & . & . & . & . & . & . & . & . & . & . \\
. & 1 & 1 & . & . & . & . & . & . & . & . & . & . & . & . \\
. & 1 & . & . & . & . & . & . & . & . & . & . & . & . & . \\
1 & . & . & . & . & . & . & . & . & . & . & . & . & . & .
\end{bmatrix},
\end{equation*}
where the entry in the $i$-th row from the top and the $j$-th column
from the left (starting from $0$) gives the coefficient of $q^i t^j$,
and dots denote zeroes.  Observe that this matrix is symmetric, a fact
that does not follow obviously from~\eqref{eq:rat-qtcat}.  We state
this symmetry as a conjecture.

\begin{conj}\label{conj:ratqt-symm}
The rational $q,t$-Catalan number is
jointly symmetric in $q$ and $t$: $$\Cat_{a,b}(q,t)=\Cat_{a,b}(t,q).$$
\end{conj} 
The weaker symmetry property $\Cat_{a,b}(q,1)=\Cat_{a,b}(1,q)$
would follow from the definition if we knew that $\sweep$ were
a bijection on $\patD(N^a E^b)$. This is known for certain choices
of $a,b$ (see~\cite{sweep}), but not for general coprime $a,b$.

\begin{example}
  The following table shows the contributions of each path $D\in
  \patD(N^3 E^5)$ to $\Cat_{3,5}(q,t)$.  Horizontal divisions
  demarcate the orbits of $\patD(N^3 E^5)$ under $\sweep$.
  \begin{table}[h]
  \centering
  \begin{tabular}{>{$}c<{$}>{$}c<{$}>{$}c<{$}>{$}c<{$}}
    D & \sweep(D) & \area(D) & \area(\sweep(D))\\\toprule
    NENENEEE & NNENEEEE & 1 & 3\\
    NNENEEEE & NENENEEE & 3 & 1\\\midrule
    NENEENEE & NNNEEEEE & 0 & 4\\
    NNNEEEEE & NENEENEE & 4 & 0\\\midrule
    NENNEEEE & NNEEENEE & 2 & 1\\
    NNEEENEE & NNEENEEE & 1 & 2\\
    NNEENEEE & NENNEEEE & 2 & 2\\\bottomrule
  \end{tabular}
  \end{table}
  We conclude that $\Cat_{3,5}(q,t) = q^4 + q^3t + q^2t^2 + q^2t + qt^2 + qt^3 + t^4$.
\end{example}

Upon substituting $t=1/q$ and rescaling by the maximum value of
$\area$ (see~\S\ref{subsec:PFab-graded}), 
we conjecture that the rational $q,t$-Catalan number
specializes to the rational $q$-Catalan number
studied in \S\ref{sec:rat-qcat}.

\begin{conj}\label{conj:qtcat-spec}
\begin{equation*}
  q^{(a-1)(b-1)/2} \Cat_{a,b}(q,1/q)=\frac{1}{[a+b]_q}\dqbin{a+b}{a,b}{q}.  
\end{equation*}
\end{conj}

\begin{remark}\label{rem:rat-qcat}
  One can use~\cite[Thm. 16]{LW-ptnid} and~\cite{mazin} along with the
  relationships among varations on the $\sweep$ map found in~\cite{sweep}
  to prove the formula $h_{\rr,\s}^+(\mu) = (\s-1)(\rr-1)/2 -
  \area(\sweep(\mu))$.  This equality implies that
  Conjectures~\ref{conj:rat-qcat} and~\ref{conj:qtcat-spec} are
  equivalent.
  Conjectures~\ref{conj:rat-qcat} and~\ref{conj:qtcat-spec} have been
  checked computationally~\cite{sage-worksheet} for $a,b\leq 12$.
\end{remark}

\begin{remark}
Rational $q,t$-Catalan numbers have been defined three times before:
by the second and third authors~\cite[\S7, Def. 21]{LW-ptnid}; 
by Gorsky and Mazin~\cite{GM-jacI,GM-jacII}; 
and by Armstrong, Hanusa, and Jones~\cite{AHJ}. 
Although these definitions use different notation and constructions,
they lead to the same $q,t$-polynomials.  
We give the details of this equivalence in~\cite{sweep}.
\end{remark}

\subsection{A Rational Frobenius Series} 
\label{subsec:RatFrob}

This section develops $q,t$-analogues of rational parking
functions (labeled rational Dyck paths).  The goal is to define a
representation of $S_a$ with a symmetric bi-grading such that the
bi-graded Hilbert series is a $q,t$-analogue of $b^{a-1}$. We give a
combinatorial formula for the Frobenius characteristic, rather than
constructing the representation itself.

To each $(a,b)$-parking function $P\in \PF_{a,b}$, we assign three
pieces of data: a $q$-weight, a $t$-weight and a subset of $[a-1]$
(used to index a fundamental quasisymmetric function).  
The $q$-weight is the usual $\area$ statistic on the
underlying Dyck path of $P$.  The $t$-weight is a ``$\dinv$''
statistic obtained as follows. 

Since $\gcd(a,b) = 1$, there exist integers $x$ and $y$ such that $xa
+ yb = 1$.  Furthermore, any integer solution $x'$ and $y'$ to
$x'a+y'b=1$ satisfies $x'=x+kb$ and $y'=y-ka$ for some integer $k$.
It follows that there is a unique such solution with $-b < x \leq 0$
and $0 \leq y < a$.  We then construct a $(|xa|,yb)$-parking function
$P'\in \PF_{|xa|,yb}$ by replacing each north step with $|x|$ copies of
itself and each east step with $y$ copies of itself.  The repeated
north steps retain their original labels from $\{1,2,\ldots,a\}$ (so
that the north-step labels now come from the multiset containing $|x|$
copies of each number from $1$ to $a$).  Since the last step of $P'$
is an east step and $|xa|=-xa = yb-1$, the first $|xa|+yb$ steps of $P'$
naturally encode a parking function $P''\in \PF_{|xa|}$.

For an $(a,b)$-parking function $P$ with underlying $(a,b)$-Dyck path $D$,
define $d(P) = \area(\sweep(D))$.  As we run over all parking
functions with a given underlying $(a,b)$-Dyck path $D$, there is a
maximum value $m(P)$ of $\dinv(P'')$.  We set $\dinv(P) = \dinv(P'') +
d(P) - m(P)$.  Note that the adjustment ensures that the maximum value
of $\dinv(P)$ is $d(P)$ for fixed underlying $(a,b)$-Dyck path $D$.

Finally, we compute the subset of $[a-1]$ assigned to $P$ as follows.
Recall from the proof of Proposition~\ref{prop:multinomial} 
that each lattice point $(x,y)$
receives a level $by-ax$.  Then the north (resp. east) steps in $P$ inherit
the level of their bottom (resp. left) endpoints.  The \textbf{diagonal
  reading word} of $P$, $\drw(P)$, is then defined to be the word obtained
by reading the labels of the north steps in increasing order according
to level.  The \textbf{inverse descent set} $\IDes(P)$ is computed
from $\drw(P)$ exactly as in \S\ref{subsec:shuffle}.

\begin{definition}\label{def:abpf}
  Define
\begin{equation*}
  \PF_{a,b}(q,t)=\sum_{P\in\PF_{a,b}} q^{\area(P)}t^{\dinv(P)}F_{n,\IDes(P)}(\x).
\end{equation*}
\end{definition}

\begin{remark}\label{rem:llt}
  One can give an equivalent definition of the rational parking
  functions in terms of the \textbf{LLT polynomials} of Lascoux,
  Leclerc, and Thibon~\cite{LLT,LT00} by assigning a suitably scaled
  LLT polynomial to each rational Dyck path.  While the LLT definition
  somewhat obscures the combinatorics, it has the advantage of
  showing immediately that the $\PF_{a,b}(q,t)$ are Schur-positive
  symmetric functions.
\end{remark}

\begin{remark}
  The diagonal reading word of a parking function $P\in \PF_{n,n+1}$
  is the reverse of the diagonal reading word obtained when $P$ is
  viewed as an element of $\PF_n$.  Our convention for $\PF_n$ is
  consistent with the existing literature.  Our convention for
  $\PF_{a,b}$ is motivated by the combinatorics.  

  Let $\omega$ denote the endomorphism on the ring of symmetric
  functions determined by $p_k \mapsto (-1)^{k-1}p_k$.  It is a
  standard fact that $\omega(s_\lambda) = s_{\lambda'}$ where
  $\lambda'$ denotes the transpose of the partition $\lambda$ (see,
  e.g.,~\cite[Theorem 1.20.2]{haglund}).  Choosing the opposite
  convention for computing the diagonal reading word for rational
  parking functions (i.e., reading the labels of the north steps in
  \emph{decreasing} order) leads to $\omega(\PF_{a,b}(q,t))$,
\end{remark}

Presumably $\PF_{a,b}(q,t)$ is the Frobenius series of some naturally
occurring bi-graded version of $\epsilon\otimes\PF_{a,b}$. We do not
know any natural construction of this bi-graded $S_a$-module analogous
to the diagonal coinvariant ring, but see Oblomkov and
Yun~\cite{OY}. For now, we make the following conjectures.

\begin{conj}\label{conj:abpf}
  Let $\Schro_{a,b;k}(q,t)$ denote the coefficient of the hook Schur
  function $s_{(k+1,1^{a-k-1})}$ in $\PF_{a,b}(q,t)$. We conjecture
  that the Schur-positive symmetric function $\PF_{a,b}(q,t)$
  satisfies the following properties.
\begin{enumerate}
\item $\PF_{a,b}(q,t)=\PF_{a,b}(t,q)$.\label{conj:ptsymm}
\item $\langle q^{(a-1)(b-1)/2}\PF_{a,b}(q,1/q),h_{(1^a)}\rangle
          =[b]_q^{a-1}$.
\item $q^{(2ak-k-k^2+ba-a^2-b+1)/2}\,\Schro_{a,b;k}(q,1/q)
    =\frac{1}{[b]_q}\dqbin{a-1}{k}{q} \dqbin{b+k}{a}{q}$.
\end{enumerate}
\end{conj}

\begin{example}
  The following table illustrates the three parking functions that
  contribute to $\PF_{2,3}(q,t)$.  Column headings refer to the
  paragraphs leading up to Definition~\ref{def:abpf}.  Note that
  north-step labels in the second column are listed from bottom to
  top.  Also note that $|x|=y=1$ (see \S\ref{subsec:RatFrob}).
\begin{table}[h]
  \centering
  \begin{tabular}{>{$}c<{$}>{$}c<{$}>{$}c<{$}cccc>{$}c<{$}>{$}c<{$}}\toprule 
  P & \text{labels} & \area & $\dinv(P'')$ & $d(P)$ & $m(P)$ & $\dinv(P)$&\IDes & \text{Contribution}\\\midrule
  NNEEE & 12 & 1 & $0$ & $0$ & $0$ & $0$ &\emptyset & qF_{2,\emptyset}(\x)\\
  NENEE & 12 & 0 & $1$ & $1$ & $1$ & $1$ &\emptyset & tF_{2,\emptyset}(\x)\\
  NENEE & 21 & 0 & $0$ & $1$ & $1$ & $0$ &\{1\} & F_{2,\{1\}}(\x)\\\bottomrule
    \end{tabular}
  \end{table}
  Since $F_{2,\emptyset} = s_2$ and $F_{2,\{1\}} = s_{1,1}$, we
  conclude that $\PF_{2,3}(q,t) = (q+t)s_2 + s_{1,1}$.
\end{example}

\begin{example}
  In this example we compute the contribution of a given parking
  function to $\PF_{5,8}(q,t)$.  Let $D = NNENNEENEEEEE$ and let $P$
  be the parking function with underlying path $D$ and with north-step
  labels $4,5,1,3,2$ when read from bottom to top.  We see that
  $\area(P) = 9$.  Also, $\sweep(D) = NENENEENNEEEE$, so $d(P) =
  \area(\sweep(D)) = 3$.  The values of $|x|$ and $y$ used in the
  construction of $P''$ are $3$ and $2$, respectively.  So $P''$ has
  underlying Dyck path $D'' = N^6E^2N^6E^4N^3E^{9}$.  A direct
  computation shows that $\dinv(P'') = 5$.  One can also check that
  $\dinv$ achieves a maximum value of $6$ on $D''$ with a label
  ordering such as $1,2,3,5,4$.  Hence $\dinv(P) = 5 + 3 - 6 = 2$.
  Finally, since $\drw(P) = 45123$, we see that $\IDes(P) = \{3\}$.
  In summary, the ultimate contribution of $P$ to $\PF_{5,8}(q,t)$ is
  $q^9t^2F_{5,\{3\}}(\x)$.
\end{example}

Part~(\ref{conj:ptsymm}) of the conjecture implies that the
coefficient of each Schur function in $\PF_{a,b}(q,t)$ is symmetric in
$q$ and $t$. For example, here is the coefficient of $s_{(3,2)}$ in
$\PF_{5,8}(q,t)$, where the entry in the $i$-th row and $j$-column
(starting from $0$) gives the coefficient of $q^it^j$ in
$\PF_{5,8}(q,t)$:
\begin{equation*}
\left[\begin{array}{rrrrrrrrrrrrr}
. & . & . & . & . & . & . & . & 1 & 1 & 1 & 1 & 1 \\
. & . & . & . & . & . & 2 & 4 & 4 & 3 & 2 & 1 & . \\
. & . & . & . & 1 & 5 & 8 & 6 & 4 & 2 & 1 & . & . \\
. & . & . & 2 & 7 & 10 & 7 & 4 & 2 & 1 & . & . & . \\
. & . & 1 & 7 & 11 & 7 & 4 & 2 & 1 & . & . & . & . \\
. & . & 5 & 10 & 7 & 4 & 2 & 1 & . & . & . & . & . \\
. & 2 & 8 & 7 & 4 & 2 & 1 & . & . & . & . & . & . \\
. & 4 & 6 & 4 & 2 & 1 & . & . & . & . & . & . & . \\
1 & 4 & 4 & 2 & 1 & . & . & . & . & . & . & . & . \\
1 & 3 & 2 & 1 & . & . & . & . & . & . & . & . & . \\
1 & 2 & 1 & . & . & . & . & . & . & . & . & . & . \\
1 & 1 & . & . & . & . & . & . & . & . & . & . & . \\
1 & . & . & . & . & . & . & . & . & . & . & . & .
\end{array}\right].
\end{equation*}

Conjecture~\ref{conj:abpf} has been checked for coprime $a,b\leq 9$.  
Note that the Schur expansion can be obtained combinatorially from the
fundamental quasisymmetric function expansion given in the definition
of rational parking functions by the methods in~\cite{elw}.
We end by posing the following question.

\begin{question}
  We know $\PF_{a,b}(q,t)$ is the Frobenius characteristic of some
  doubly-graded $S_a$-module of dimension $b^{a-1}$, whereas
  $\PF_{b,a}(q,t)$ is the Frobenius characteristic of some
  doubly-graded $S_b$-module of dimension $a^{b-1}$.  
  What is the algebraic and combinatorial relationship between
  $\PF_{a,b}(q,t)$ and $\PF_{b,a}(q,t)$?
\end{question}

\section{Appendix: Computations}
\label{sec:comp}

The following expansions were computed using Sage~\cite{sage}.  The
Sage worksheet we used to do this can be found on the third author's web
page~\cite{sage-worksheet}.


\setlength{\arraycolsep}{3pt}

\subsection{Rational Catalan numbers}
\label{subsec:cat-comp}
\[
\Cat_{2,3}(q,t) = \Cat_{3,2}(q,t) = 
\begin{bmatrix}
 . &  1 \\
 1 &  . \\
\end{bmatrix},
\]

\[
\Cat_{3,5}(q,t) = \Cat_{5,3}(q,t) = 
\begin{bmatrix}
 . &  . &  . &  . &  1 \\
 . &  . &  1 &  1 &  . \\
 . &  1 &  1 &  . &  . \\
 . &  1 &  . &  . &  . \\
 1 &  . &  . &  . &  . \\
\end{bmatrix},
\]

\[
\Cat_{3,7}(q,t) = \Cat_{7,3}(q,t) = 
\begin{bmatrix}
 . &  . &  . &  . &  . &  . &  1 \\
 . &  . &  . &  . &  1 &  1 &  . \\
 . &  . &  1 &  1 &  1 &  . &  . \\
 . &  . &  1 &  1 &  . &  . &  . \\
 . &  1 &  1 &  . &  . &  . &  . \\
 . &  1 &  . &  . &  . &  . &  . \\
 1 &  . &  . &  . &  . &  . &  . \\
\end{bmatrix},
\]

\[
\Cat_{4,7}(q,t) = \Cat_{7,4}(q,t) =
\begin{bmatrix}
 . &  . &  . &  . &  . &  . &  . &  . &  . &  1 \\
 . &  . &  . &  . &  . &  . &  1 &  1 &  1 &  . \\
 . &  . &  . &  . &  1 &  2 &  1 &  1 &  . &  . \\
 . &  . &  . &  1 &  2 &  1 &  1 &  . &  . &  . \\
 . &  . &  1 &  2 &  1 &  1 &  . &  . &  . &  . \\
 . &  . &  2 &  1 &  1 &  . &  . &  . &  . &  . \\
 . &  1 &  1 &  1 &  . &  . &  . &  . &  . &  . \\
 . &  1 &  1 &  . &  . &  . &  . &  . &  . &  . \\
 . &  1 &  . &  . &  . &  . &  . &  . &  . &  . \\
 1 &  . &  . &  . &  . &  . &  . &  . &  . &  . \\
\end{bmatrix},
\]

\[
\Cat_{5,8}(q,t) = \Cat_{8,5}(q,t) = 
\begin{bmatrix}
 . &  . &  . &  . &  . &  . &  . &  . &  . &  . &  . &  . &  . &  . &  1 \\
 . &  . &  . &  . &  . &  . &  . &  . &  . &  . &  1 &  1 &  1 &  1 &  . \\
 . &  . &  . &  . &  . &  . &  . &  . &  2 &  2 &  2 &  1 &  1 &  . &  . \\
 . &  . &  . &  . &  . &  . &  1 &  3 &  3 &  2 &  1 &  1 &  . &  . &  . \\
 . &  . &  . &  . &  . &  2 &  4 &  3 &  2 &  1 &  1 &  . &  . &  . &  . \\
 . &  . &  . &  . &  2 &  4 &  3 &  2 &  1 &  1 &  . &  . &  . &  . &  . \\
 . &  . &  . &  1 &  4 &  3 &  2 &  1 &  1 &  . &  . &  . &  . &  . &  . \\
 . &  . &  . &  3 &  3 &  2 &  1 &  1 &  . &  . &  . &  . &  . &  . &  . \\
 . &  . &  2 &  3 &  2 &  1 &  1 &  . &  . &  . &  . &  . &  . &  . &  . \\
 . &  . &  2 &  2 &  1 &  1 &  . &  . &  . &  . &  . &  . &  . &  . &  . \\
 . &  1 &  2 &  1 &  1 &  . &  . &  . &  . &  . &  . &  . &  . &  . &  . \\
 . &  1 &  1 &  1 &  . &  . &  . &  . &  . &  . &  . &  . &  . &  . &  . \\
 . &  1 &  1 &  . &  . &  . &  . &  . &  . &  . &  . &  . &  . &  . &  . \\
 . &  1 &  . &  . &  . &  . &  . &  . &  . &  . &  . &  . &  . &  . &  . \\
 1 &  . &  . &  . &  . &  . &  . &  . &  . &  . &  . &  . &  . &  . &  . \\
\end{bmatrix}.
\]

\subsection{Rational Parking Functions}
\label{subsec:rpf-comp}
\begin{gather*}
\PF_{2,3}(q,t) = 
\begin{bmatrix}
. & 1 \\
1 & . \\
\end{bmatrix}
s_{(2)} + \
\begin{bmatrix}
1 \\
\end{bmatrix}
s_{(1, 1)},
\end{gather*}

\begin{gather*}
\PF_{2,5}(q,t) = 
\begin{bmatrix}
. & . & 1 \\
. & 1 & . \\
1 & . & . \\
\end{bmatrix}
s_{(2)} + \
\begin{bmatrix}
. & 1 \\
1 & . \\
\end{bmatrix}
s_{(1, 1)},
\end{gather*}

\begin{gather*}
\PF_{3,5}(q,t) = 
\begin{bmatrix}
. & . & . & . & 1 \\
. & . & 1 & 1 & . \\
. & 1 & 1 & . & . \\
. & 1 & . & . & . \\
1 & . & . & . & . \\
\end{bmatrix}
s_{(3)} + \
\begin{bmatrix}
. & . & 1 & 1 \\
. & 2 & 1 & . \\
1 & 1 & . & . \\
1 & . & . & . \\
\end{bmatrix}
s_{(2, 1)} + \
\begin{bmatrix}
. & 1 \\
1 & . \\
\end{bmatrix}
s_{(1, 1, 1)},
\end{gather*}

\begin{gather*}
\PF_{5,3}(q,t) = 
\begin{bmatrix}
. & . & . & . & 1 \\
. & . & 1 & 1 & . \\
. & 1 & 1 & . & . \\
. & 1 & . & . & . \\
1 & . & . & . & . \\
\end{bmatrix}
s_{(5)} + \
\begin{bmatrix}
. & . & 1 & 1 \\
. & 2 & 1 & . \\
1 & 1 & . & . \\
1 & . & . & . \\
\end{bmatrix}
s_{(4, 1)} + \
\begin{bmatrix}
. & 1 & 1 \\
1 & 1 & . \\
1 & . & . \\
\end{bmatrix}
s_{(3, 2)} + \
\begin{bmatrix}
. & 1 \\
1 & . \\
\end{bmatrix}
s_{(3, 1, 1)} + \
\begin{bmatrix}
1 \\
\end{bmatrix}
s_{(2, 2, 1)},
\end{gather*}

\begin{gather*}
\PF_{4,7}(q,t) = 
\begin{bmatrix}
. & . & . & . & . & . & . & . & . & 1 \\
. & . & . & . & . & . & 1 & 1 & 1 & . \\
. & . & . & . & 1 & 2 & 1 & 1 & . & . \\
. & . & . & 1 & 2 & 1 & 1 & . & . & . \\
. & . & 1 & 2 & 1 & 1 & . & . & . & . \\
. & . & 2 & 1 & 1 & . & . & . & . & . \\
. & 1 & 1 & 1 & . & . & . & . & . & . \\
. & 1 & 1 & . & . & . & . & . & . & . \\
. & 1 & . & . & . & . & . & . & . & . \\
1 & . & . & . & . & . & . & . & . & . \\
\end{bmatrix}
s_{(4)} + \
\begin{bmatrix}
. & . & . & . & . & . & 1 & 1 & 1 \\
. & . & . & . & 2 & 3 & 2 & 1 & . \\
. & . & 1 & 3 & 4 & 2 & 1 & . & . \\
. & . & 3 & 4 & 2 & 1 & . & . & . \\
. & 2 & 4 & 2 & 1 & . & . & . & . \\
. & 3 & 2 & 1 & . & . & . & . & . \\
1 & 2 & 1 & . & . & . & . & . & . \\
1 & 1 & . & . & . & . & . & . & . \\
1 & . & . & . & . & . & . & . & . \\
\end{bmatrix}
s_{(3, 1)} + \\ \
\begin{bmatrix}
. & . & . & . & . & 1 & . & 1 \\
. & . & . & 1 & 2 & 1 & 1 & . \\
. & . & 1 & 3 & 1 & 1 & . & . \\
. & 1 & 3 & 1 & 1 & . & . & . \\
. & 2 & 1 & 1 & . & . & . & . \\
1 & 1 & 1 & . & . & . & . & . \\
. & 1 & . & . & . & . & . & . \\
1 & . & . & . & . & . & . & . \\
\end{bmatrix}
s_{(2, 2)} + \
\begin{bmatrix}
. & . & . & . & 1 & 1 & 1 \\
. & . & 1 & 3 & 2 & 1 & . \\
. & 1 & 3 & 2 & 1 & . & . \\
. & 3 & 2 & 1 & . & . & . \\
1 & 2 & 1 & . & . & . & . \\
1 & 1 & . & . & . & . & . \\
1 & . & . & . & . & . & . \\
\end{bmatrix}
s_{(2, 1, 1)} + \
\begin{bmatrix}
. & . & . & 1 \\
. & 1 & 1 & . \\
. & 1 & . & . \\
1 & . & . & . \\
\end{bmatrix}
s_{(1, 1, 1, 1)},
\end{gather*}

\begin{gather*}
\PF_{7,4}(q,t) = 
\begin{bmatrix}
. & . & . & . & . & . & . & . & . & 1 \\
. & . & . & . & . & . & 1 & 1 & 1 & . \\
. & . & . & . & 1 & 2 & 1 & 1 & . & . \\
. & . & . & 1 & 2 & 1 & 1 & . & . & . \\
. & . & 1 & 2 & 1 & 1 & . & . & . & . \\
. & . & 2 & 1 & 1 & . & . & . & . & . \\
. & 1 & 1 & 1 & . & . & . & . & . & . \\
. & 1 & 1 & . & . & . & . & . & . & . \\
. & 1 & . & . & . & . & . & . & . & . \\
1 & . & . & . & . & . & . & . & . & . \\
\end{bmatrix}
s_{(7)} + \
\begin{bmatrix}
. & . & . & . & . & . & 1 & 1 & 1 \\
. & . & . & . & 2 & 3 & 2 & 1 & . \\
. & . & 1 & 3 & 4 & 2 & 1 & . & . \\
. & . & 3 & 4 & 2 & 1 & . & . & . \\
. & 2 & 4 & 2 & 1 & . & . & . & . \\
. & 3 & 2 & 1 & . & . & . & . & . \\
1 & 2 & 1 & . & . & . & . & . & . \\
1 & 1 & . & . & . & . & . & . & . \\
1 & . & . & . & . & . & . & . & . \\
\end{bmatrix}
s_{(6, 1)} + \\ \
\begin{bmatrix}
. & . & . & . & 1 & 2 & 1 & 1 \\
. & . & 1 & 3 & 4 & 2 & 1 & . \\
. & 1 & 4 & 5 & 2 & 1 & . & . \\
. & 3 & 5 & 2 & 1 & . & . & . \\
1 & 4 & 2 & 1 & . & . & . & . \\
2 & 2 & 1 & . & . & . & . & . \\
1 & 1 & . & . & . & . & . & . \\
1 & . & . & . & . & . & . & . \\
\end{bmatrix}
s_{(5, 2)} + \
\begin{bmatrix}
. & . & . & . & 1 & 1 & 1 \\
. & . & 1 & 3 & 2 & 1 & . \\
. & 1 & 3 & 2 & 1 & . & . \\
. & 3 & 2 & 1 & . & . & . \\
1 & 2 & 1 & . & . & . & . \\
1 & 1 & . & . & . & . & . \\
1 & . & . & . & . & . & . \\
\end{bmatrix}
s_{(5, 1, 1)} + \
\begin{bmatrix}
. & . & . & 1 & 1 & 1 & 1 \\
. & . & 2 & 3 & 2 & 1 & . \\
. & 2 & 4 & 2 & 1 & . & . \\
1 & 3 & 2 & 1 & . & . & . \\
1 & 2 & 1 & . & . & . & . \\
1 & 1 & . & . & . & . & . \\
1 & . & . & . & . & . & . \\
\end{bmatrix}
s_{(4, 3)} + \\ \
\begin{bmatrix}
. & . & 1 & 2 & 2 & 1 \\
. & 2 & 4 & 3 & 1 & . \\
1 & 4 & 3 & 1 & . & . \\
2 & 3 & 1 & . & . & . \\
2 & 1 & . & . & . & . \\
1 & . & . & . & . & . \\
\end{bmatrix}
s_{(4, 2, 1)} + \
\begin{bmatrix}
. & . & . & 1 \\
. & 1 & 1 & . \\
. & 1 & . & . \\
1 & . & . & . \\
\end{bmatrix}
s_{(4, 1, 1, 1)} + \
\begin{bmatrix}
. & . & 1 & 1 & 1 \\
. & 2 & 2 & 1 & . \\
1 & 2 & 1 & . & . \\
1 & 1 & . & . & . \\
1 & . & . & . & . \\
\end{bmatrix}
s_{(3, 3, 1)} + \\ \
\begin{bmatrix}
. & 1 & 1 & 1 \\
1 & 1 & 1 & . \\
1 & 1 & . & . \\
1 & . & . & . \\
\end{bmatrix}
s_{(3, 2, 2)} + \
\begin{bmatrix}
. & 1 & 1 \\
1 & 1 & . \\
1 & . & . \\
\end{bmatrix}
s_{(3, 2, 1, 1)} + \
\begin{bmatrix}
1 \\
\end{bmatrix}
s_{(2, 2, 2, 1)},
\end{gather*}

\begin{gather*}
\PF_{5,8}(q,t) = 
\begin{bmatrix}
. & . & . & . & . & . & . & . & . & . & . & . & . & . & 1 \\
. & . & . & . & . & . & . & . & . & . & 1 & 1 & 1 & 1 & . \\
. & . & . & . & . & . & . & . & 2 & 2 & 2 & 1 & 1 & . & . \\
. & . & . & . & . & . & 1 & 3 & 3 & 2 & 1 & 1 & . & . & . \\
. & . & . & . & . & 2 & 4 & 3 & 2 & 1 & 1 & . & . & . & . \\
. & . & . & . & 2 & 4 & 3 & 2 & 1 & 1 & . & . & . & . & . \\
. & . & . & 1 & 4 & 3 & 2 & 1 & 1 & . & . & . & . & . & . \\
. & . & . & 3 & 3 & 2 & 1 & 1 & . & . & . & . & . & . & . \\
. & . & 2 & 3 & 2 & 1 & 1 & . & . & . & . & . & . & . & . \\
. & . & 2 & 2 & 1 & 1 & . & . & . & . & . & . & . & . & . \\
. & 1 & 2 & 1 & 1 & . & . & . & . & . & . & . & . & . & . \\
. & 1 & 1 & 1 & . & . & . & . & . & . & . & . & . & . & . \\
. & 1 & 1 & . & . & . & . & . & . & . & . & . & . & . & . \\
. & 1 & . & . & . & . & . & . & . & . & . & . & . & . & . \\
1 & . & . & . & . & . & . & . & . & . & . & . & . & . & . \\
\end{bmatrix}
s_{(5)} + \
\begin{bmatrix}
. & . & . & . & . & . & . & . & . & . & 1 & 1 & 1 & 1 \\
. & . & . & . & . & . & . & 1 & 3 & 4 & 3 & 2 & 1 & . \\
. & . & . & . & . & 1 & 4 & 7 & 6 & 4 & 2 & 1 & . & . \\
. & . & . & . & 2 & 7 & 9 & 7 & 4 & 2 & 1 & . & . & . \\
. & . & . & 2 & 8 & 10 & 7 & 4 & 2 & 1 & . & . & . & . \\
. & . & 1 & 7 & 10 & 7 & 4 & 2 & 1 & . & . & . & . & . \\
. & . & 4 & 9 & 7 & 4 & 2 & 1 & . & . & . & . & . & . \\
. & 1 & 7 & 7 & 4 & 2 & 1 & . & . & . & . & . & . & . \\
. & 3 & 6 & 4 & 2 & 1 & . & . & . & . & . & . & . & . \\
. & 4 & 4 & 2 & 1 & . & . & . & . & . & . & . & . & . \\
1 & 3 & 2 & 1 & . & . & . & . & . & . & . & . & . & . \\
1 & 2 & 1 & . & . & . & . & . & . & . & . & . & . & . \\
1 & 1 & . & . & . & . & . & . & . & . & . & . & . & . \\
1 & . & . & . & . & . & . & . & . & . & . & . & . & . \\
\end{bmatrix}
s_{(4, 1)} + \\ \
\begin{bmatrix}
. & . & . & . & . & . & . & . & 1 & 1 & 1 & 1 & 1 \\
. & . & . & . & . & . & 2 & 4 & 4 & 3 & 2 & 1 & . \\
. & . & . & . & 1 & 5 & 8 & 6 & 4 & 2 & 1 & . & . \\
. & . & . & 2 & 7 & 10 & 7 & 4 & 2 & 1 & . & . & . \\
. & . & 1 & 7 & 11 & 7 & 4 & 2 & 1 & . & . & . & . \\
. & . & 5 & 10 & 7 & 4 & 2 & 1 & . & . & . & . & . \\
. & 2 & 8 & 7 & 4 & 2 & 1 & . & . & . & . & . & . \\
. & 4 & 6 & 4 & 2 & 1 & . & . & . & . & . & . & . \\
1 & 4 & 4 & 2 & 1 & . & . & . & . & . & . & . & . \\
1 & 3 & 2 & 1 & . & . & . & . & . & . & . & . & . \\
1 & 2 & 1 & . & . & . & . & . & . & . & . & . & . \\
1 & 1 & . & . & . & . & . & . & . & . & . & . & . \\
1 & . & . & . & . & . & . & . & . & . & . & . & . \\
\end{bmatrix}
s_{(3, 2)} + \
\begin{bmatrix}
. & . & . & . & . & . & . & 1 & 1 & 2 & 1 & 1 \\
. & . & . & . & . & 2 & 4 & 5 & 4 & 2 & 1 & . \\
. & . & . & 1 & 4 & 8 & 7 & 5 & 2 & 1 & . & . \\
. & . & 1 & 5 & 10 & 8 & 5 & 2 & 1 & . & . & . \\
. & . & 4 & 10 & 8 & 5 & 2 & 1 & . & . & . & . \\
. & 2 & 8 & 8 & 5 & 2 & 1 & . & . & . & . & . \\
. & 4 & 7 & 5 & 2 & 1 & . & . & . & . & . & . \\
1 & 5 & 5 & 2 & 1 & . & . & . & . & . & . & . \\
1 & 4 & 2 & 1 & . & . & . & . & . & . & . & . \\
2 & 2 & 1 & . & . & . & . & . & . & . & . & . \\
1 & 1 & . & . & . & . & . & . & . & . & . & . \\
1 & . & . & . & . & . & . & . & . & . & . & . \\
\end{bmatrix}
s_{(3, 1, 1)} + \\ \
\begin{bmatrix}
. & . & . & . & . & . & 1 & 1 & 1 & 1 & 1 \\
. & . & . & . & 1 & 3 & 4 & 3 & 2 & 1 & . \\
. & . & . & 2 & 6 & 6 & 4 & 2 & 1 & . & . \\
. & . & 2 & 7 & 7 & 4 & 2 & 1 & . & . & . \\
. & 1 & 6 & 7 & 4 & 2 & 1 & . & . & . & . \\
. & 3 & 6 & 4 & 2 & 1 & . & . & . & . & . \\
1 & 4 & 4 & 2 & 1 & . & . & . & . & . & . \\
1 & 3 & 2 & 1 & . & . & . & . & . & . & . \\
1 & 2 & 1 & . & . & . & . & . & . & . & . \\
1 & 1 & . & . & . & . & . & . & . & . & . \\
1 & . & . & . & . & . & . & . & . & . & . \\
\end{bmatrix}
s_{(2, 2, 1)} + \
\begin{bmatrix}
. & . & . & . & . & 1 & 1 & 1 & 1 \\
. & . & . & 1 & 3 & 3 & 2 & 1 & . \\
. & . & 2 & 4 & 4 & 2 & 1 & . & . \\
. & 1 & 4 & 4 & 2 & 1 & . & . & . \\
. & 3 & 4 & 2 & 1 & . & . & . & . \\
1 & 3 & 2 & 1 & . & . & . & . & . \\
1 & 2 & 1 & . & . & . & . & . & . \\
1 & 1 & . & . & . & . & . & . & . \\
1 & . & . & . & . & . & . & . & . \\
\end{bmatrix}
s_{(2, 1, 1, 1)} + \
\begin{bmatrix}
. & . & . & . & 1 \\
. & . & 1 & 1 & . \\
. & 1 & 1 & . & . \\
. & 1 & . & . & . \\
1 & . & . & . & . \\
\end{bmatrix}
s_{(1, 1, 1, 1, 1)}.
\end{gather*}

\section*{Acknowledgments}

The authors gratefully acknowledge discussions with Eugene Gorsky, Jim
Haglund, Mark Haiman, Mikhail Mazin and Michelle Wachs.


\begin{thebibliography}{99}
\bibitem{armstrong} Drew Armstrong, ``Hyperplane arrangements and diagonal harmonics,'' \emph{J. Comb.}, \textbf{4} (2013), no. 2, 157--190.

\bibitem{AHJ} Drew Armstrong, Christopher R. H. Hanusa and Brant
  C. Jones, ``Results and conjectures on simultaneous core
  partitions,'' \texttt{http://arxiv.org/abs/1308.0572}.

\bibitem{ARW} Drew Armstrong, Brendon Rhoades, and Nathan Williams, ``Rational associahedra and noncrossing partitions,'' \texttt{http://arxiv.org/abs/1305.7286}.

\bibitem{sweep} Drew Armstrong, Nicholas A. Loehr, and Gregory
  S. Warrington, ``Sweep maps: A continuous family of sorting
  algorithms,'' in preparation.

\bibitem{bizley} M.T.L. Bizley, ``Derivation of a new formula for the number of minimal lattice paths from $(0,0)$ to $(km,kn)$ having just $t$ contacts with the line $my=nx$ and having no points above this line; and a proof of Grossman's formula for the number of paths which may touch but do not rise above this line,'' \emph{J. Inst. Actuar.} \textbf{80} (1954), 55--62.

\bibitem{borel} Armand Borel, ``Sur la cohomologie des espaces
  fibr\'es principaux et des espaces homog\`enes de groupes de Lie
  compacts,'' 
\emph{Ann. of Math.} \textbf{57} no. 2 (1953), 115--207.

\bibitem{cayley}
A.~Cayley, ``A theorem on trees,'' Quart. J. Math. {\bf 23} (1889), 376--378.

\bibitem{chevalley} Claude Chevalley, ``Invariants of finite groups
  generated by reflections,'' \emph{Amer. J. Math.}, \textbf{77}
  (1955), 778--782.

\bibitem{elw} Eric Egge, Nicholas Loehr and Gregory S. Warrington,
  ``From quasisymmetric expansions to Schur expansions via a modified
  inverse Kostka matrix,'' European J. Combin. 31 (2010), no. 8,
  2014--2027.

\bibitem{qtcat-proof} Jim Haglund and Adriano Garsia, ``A proof of the
  $q,t$-Catalan positivity conjecture,'' 
  \emph{Discrete Math.} \textbf{256} (2002), 
  677--717.

\bibitem{GH-qtcat} Adriano Garsia and Mark Haiman, ``A remarkable $q,t$-Catalan
 sequence and $q$-Lagrange inversion,'' \emph{J. Algebraic Combin.}
 \textbf{5} (1996), 191--244.

\bibitem{GM-jacI} Evgeny Gorsky and Mikhail Mazin,
 ``Compactified Jacobians and $q,t$-Catalan numbers, I,''
 \emph{J. Combin. Theory Ser. A} \textbf{120} (2013), 49--63.
 Also online at \texttt{arXiv:1105.1151v2}.

\bibitem{GM-jacII} Evgeny Gorsky and Mikhail Mazin,
 ``Compactified Jacobians and $q,t$-Catalan numbers, II,'' 
 preprint, \texttt{arXiv:1204.5448v1}.

\bibitem{GMV} Eugene Gorsky, Mikhail Mazin, and Monica Vazirani, ``Affine permutations and rational slope parking functions,'' \texttt{http://arxiv.org/abs/1403.0303}.

\bibitem{gorsky-negut} Eugene Gorsky and Andrew Negut, ``Refined knot invariants and Hilbert schemes,'' \texttt{http://arxiv.org/abs/1304.3328v2}.

\bibitem{Groj} Ian Grojnowski and Mark Haiman, ``Affine Hecke
  algebras and positivity of LLT and Macdonald polynomials,'' preprint.

\bibitem{haglund} James Haglund, ``The $q,t$-Catalan Numbers and the Space of Diagonal Harmonics,'' University Lecture Series {\bf 41},  American Mathematical Society, Providence, RI, 2008.

\bibitem{Hag-bounce} James Haglund, ``Conjectured statistics for the
 $q,t$-Catalan numbers,'' \emph{Adv. in Math.} \textbf{175} (2003),
 319--334.

\bibitem{HHLRU} James Haglund, Mark Haiman, Nicholas Loehr,
 Jeff Remmel, and A. Ulyanov, ``A combinatorial formula for the
 character of the diagonal coinvariants,'' \emph{Duke Math. J.}
 \textbf{126} (2005), 195--232.

\bibitem{HL-park} James Haglund and Nicholas A. Loehr,
 ``A conjectured combinatorial formula for the Hilbert series for
 diagonal harmonics,'' \emph{Discrete Math.} \textbf{298} (2005), 189--204.

\bibitem{haiman-conjectures} Mark Haiman, ``Conjectures on the quotient ring by diagonal invariants,'' \emph{J. Algebraic Combin.} \textbf{3} (1994), 17--76.

\bibitem{haiman-vanish} Mark Haiman, ``Vanishing theorems and character
 formulas for the Hilbert scheme of points in the plane,'' \emph{Invent.
 Math.} \textbf{149} (2002), 371--407.

\bibitem{hikita} Tatsuyuki Hikita, ``Affine Springer fibers of type $A$ and combinatorics of diagonal coinvariants,'' \texttt{http://arxiv.org/abs/1203.5878v1}.

\bibitem{KT02} Masaki Kashiwara and Toshiyuki Tanisaki, ``Parabolic
  Kazhdan-Lusztig polynomials and Schubert varieties,''
  \emph{J. Algebra}, \textbf{249} (2002), no. 2, 306--325.

\bibitem{konheim-weiss}
A. G. Konheim and B. Weiss, ``An occupancy discipline and applications,''
\emph{SIAM J. Applied Math.} \textbf{14} (1966), 17--76.

\bibitem{kreweras}
G.~Kreweras, ``Sur les partitions non crois\'ees d'un cycle,''
\emph{Discrete Math.} {\bf 1} (1972), 333--350.

\bibitem{LLT} Alain Lascoux, Bernard Leclerc and Jean-Yves Thibon,
  ``Ribbon tableaux, Hall-Littlewood functions, quantum affine
  algebras, and unipotent varieties,'' \emph{J. Math. Phys.},
  \textbf{38} (1997), no. 2, 1041--1068.

\bibitem{LT00} Bernard Leclerc and Jean-Yves Thibon,
  ``Littlewood-Richardson coefficients and Kazhdan-Lusztig
  polynomials,'' Combinatorial methods in representation theory
  (Kyoto, 1998), \emph{Adv. Stud. Pure Math.} \textbf{28} (2000), 155--220.

\bibitem{loehr-bij} Nicholas A. Loehr, \emph{Bijective Combinatorics},
 Chapman and Hall/CRC Press (2011).

\bibitem{loehr-mcat} Nicholas A. Loehr, ``Conjectured statistics for
 the higher $q,t$-Catalan sequences,'' \emph{Electron. J. Combin.} \textbf{12}
 (2005) research paper R9; 54 pages (electronic).

\bibitem{LW-ptnid} Nicholas A. Loehr and Gregory S. Warrington,
 ``A continuous family of partition statistics equidistributed with length,''
 \emph{J. Combin. Theory Ser. A} \textbf{116} (2009), 379--403.

\bibitem{Macd} Ian Macdonald, \emph{Symmetric Functions and Hall Polynomials}
(2nd ed.), Oxford University Press (1995).

\bibitem{macmahon}
 Percy MacMahon, \emph{Combinatory Analysis},
 Cambridge University Press (1918).
 Reprinted by Chelsea, New York, 3rd edition (1984).

\bibitem{mazin} Mikhail Mazin, ``A bijective proof of
  Loehr-Warrington's formulas for the statistics
  $\mbox{ctot}_{\frac{q}{p}}$ and $\mbox{mid}_{\frac{q}{p}}$,'' preprint,
  \texttt{arXiv:1301.7452}.

\bibitem{OY} Alexei Oblomkov and Zhiwei Yun, ``Geometric
  representations of graded and rational Cherednik algebras,'' preprint.

\bibitem{pak-postnikov} Igor Pak and Alexander Postnikov, ``Enumeration of trees and one amazing representation of $S_n$,'' \emph{Proc. FPSAC '96 Conf.}, Minneapolis, MN.

\bibitem{riordan} John Riordan, ``Ballots and trees,'' \emph{J. Combinatorial Theory} \textbf{6} (1969), 408--411.

\bibitem{sagan} Bruce Sagan, \emph{The Symmetric Group:
 Representations, Combinatorial Algorithms, and Symmetric Functions}
 (second ed.), Springer, New York (2001).

\bibitem{stanley-ec2} Richard P. Stanley, \emph{Enumerative Combinatorics 
Volume 2}, Cambridge, 1999.

\bibitem{stanley-parking} Richard P. Stanley, ``Parking functions and noncrossing partitions,'' \emph{Electron. J. Combin.} \textbf{4} (1997), no. 2, Research Paper 20, approx. 14 pp.

\bibitem{sage} W. A. Stein et al., \emph{Sage Mathematics Software
    (Version 5.5)}, The Sage Development Team, 2012.
  \url{http://www.sagemath.org}.

\bibitem{sage-worksheet} Gregory S. Warrington, \emph{Sage worksheet for
    rational Catalan and rational parking function conjectures
    (available online)}, 2014.
  \url{http://www.cems.uvm.edu/~gswarrin/research/RPF.sws}.  

\bibitem{weyl} Hermann Weyl, \emph{The Classical Groups, Their
    Invariants and Representations}, 2nd. ed., Princeton Univ. Press,
  Princeton, N.J., 1946.

\end{thebibliography}
\end{document}